\documentclass[10pt,reqno]{smfart}
\usepackage{amsmath,amssymb,smfthm,stmaryrd,epic}
\usepackage[frenchb]{babel}
\usepackage[latin1]{inputenc}
\usepackage[all]{xy}
\usepackage{a4wide}

\usepackage{xypic}
\usepackage[active]{srcltx}
\usepackage{rotating}

\NumberTheoremsIn{subsection}\SwapTheoremNumbers

\begin{document}

\PUSH{macros2.tex}%
\newtheorem{prop-defi}[smfthm]{Proposition-Définition}
\newtheorem{theo-defi}[smfthm]{Théorème-définition}
\newtheorem{lem-defi}[smfthm]{Lemme-définition}
\newtheorem{notas}[smfthm]{Notations}
\newtheorem{nota}[smfthm]{Notation}
\newtheorem{defis}[smfthm]{Définitions}
\newtheorem{remas}[smfthm]{Remarques}

\newtheorem{theob}{Théorème}[section]
\def\thetheob{\arabic{section}.\arabic{theob}}
\newtheorem{propb}[theob]{Proposition}
\newtheorem{lemb}[theob]{Lemme}
\newtheorem{corob}[theob]{Corollaire}
\newtheorem{defib}[theob]{Définition}
\newtheorem{defisb}[theob]{Définitions}
\newtheorem{remab}[theob]{Remarque}

\renewcommand{\theequation}{\Roman{part}.\arabic{section}.\arabic{subsection}.\arabic{smfthm}}

\def\Am{{\mathbb A}}
\def\Fm{{\mathbb F}}
\def\Mm{{\mathbb M}}
\def\Nm{{\mathbb N}}
\def\Pm{{\mathbb P}}
\def\Qm{{\mathbb Q}}
\def\Zm{{\mathbb Z}}
\def\Dm{{\mathbb D}}
\def\Cm{{\mathbb C}}
\def\Rm{{\mathbb R}}
\def\Gm{{\mathbb G}}
\def\Lm{{\mathbb L}}
\def \Km{{\mathbb K}}

\def\ZC{{\mathcal Z}}
\def\AC{{\mathcal A}}
\def\CC{{\mathcal C}}
\def\DC{{\mathcal D}}
\def\EC{{\mathcal E}}
\def\FC{{\mathcal F}}
\def\GC{{\mathcal G}}
\def\HC{{\mathcal H}}
\def\IC{{\mathcal I}}
\def\JC{{\mathcal J}}
\def\KC{{\mathcal K}}
\def\LC{{\mathcal L}}
\def\MC{{\mathcal M}}
\def\NC{{\mathcal N}}
\def\OC{{\mathcal O}}
\def\PC{{\mathcal P}}
\def\UC{{\mathcal U}}
\def\VC{{\mathcal V}}
\def\XC{{\mathcal X}}
\def\YC{{\mathcal Y}}

\def\BF{{\mathfrak B}}
\def\AF{{\mathfrak A}}
\def\GF{{\mathfrak G}}
\def\EF{{\mathfrak E}}
\def\CF{{\mathfrak C}}
\def\DF{{\mathfrak D}}
\def\JF{{\mathfrak J}}
\def\LF{{\mathfrak L}}
\def\MF{{\mathfrak M}}
\def\NF{{\mathfrak N}}
\def\XF{{\mathfrak X}}
\def\UF{{\mathfrak U}}
\def\KF{{\mathfrak K}}

\def \longmapright#1{\smash{\mathop{\longrightarrow}\limits^{#1}}}
\def \mapright#1{\smash{\mathop{\rightarrow}\limits^{#1}}}
\def \lexp#1#2{\kern \scriptspace \vphantom{#2}^{#1}\kern-\scriptspace#2}
\def \linf#1#2{\kern \scriptspace \vphantom{#2}_{#1}\kern-\scriptspace#2}
\def \linexp#1#2#3 {\kern \scriptspace{#3}_{#1}^{#2} \kern-\scriptspace #3}

\def \a{\alpha}
\def \b{\beta}
\def \d{\delta}
\def \e{\epsilon}
\def \g{\gamma}
\def \k{\kappa}
\def \l{\lambda}
\def \m{\mu}
\def \n{\nu}
\def \o{\omega}
\def \r{\rho}
\def \s{\sigma}
\def \t{\tau}
\def \th{\theta}
\def \u {\upsilon}
\def \x{\chi}
\def \vphi {\varphi}

\let \leq=\leqslant
\let \geq=\geqslant
\def \lefto{\longleftarrow}
\def \fin{\hfill $\square$}
\let \DS=\displaystyle
\let \SS=\scriptstyle
\let \longto=\longrightarrow
\let \oo=\infty

\def \sd{\mathop{\mathrm{Sd}}\nolimits}
\def \irr{\mathop{\mathrm{Irr}}\nolimits}
\def \FH{\mathop{\mathrm{FH}}\nolimits}
\def \FPH{\mathop{\mathrm{FPH}}\nolimits}
\def \coh{\mathop{\mathrm{Coh}}\nolimits}
\def \res{\mathop{\mathrm{res}}\nolimits}
\def \op{\mathop{\mathrm{op}}\nolimits}
\def \rec {\mathop{\mathrm{rec}}\nolimits}
\def \art{\mathop{\mathrm{Art}}\nolimits}
\def \hyp {\mathop{\mathrm{Hyp}}\nolimits}
\def \cusp {\mathop{\mathrm{Cusp}}\nolimits}
\def \Iw {\mathop{\mathrm{Iw}}\nolimits}
\def \JL {\mathop{\mathrm{JL}}\nolimits}
\def \speh {\mathop{\mathrm{Speh}}\nolimits}
\def \isom {\mathop{\mathrm{Isom}}\nolimits}
\def \Vect {\mathop{\mathrm{Vect}}\nolimits}
\def \groth {\mathop{\mathrm{Groth}}\nolimits}
\def \lef {\mathop{\mathrm{Lef}}\nolimits}
\def \fix {\mathop{\mathrm{Fix}}\nolimits}
\def \hom {\mathop{\mathrm{Hom}}\nolimits}
\def \deg {\mathop{\mathrm{deg}}\nolimits}
\def \val {\mathop{\mathrm{val}}\nolimits}
\def \det {\mathop{\mathrm{det}}\nolimits}
\def \rep {\mathop{\mathrm{Rep}}\nolimits}
\def \spec {\mathop{\mathrm{Spec}}\nolimits}
\def \fr {\mathop{\mathrm{Fr}}\nolimits}
\def \frob {\mathop{\mathrm{Frob}}\nolimits}
\def \ker {\mathop{\mathrm{Ker}}\nolimits}
\def \im {\mathop{\mathrm{Im}}\nolimits}
\def \Red {\mathop{\mathrm{Red}}\nolimits}
\def \red {\mathop{\mathrm{red}}\nolimits}
\def \aut {\mathop{\mathrm{Aut}}\nolimits}
\def \diag {\mathop{\mathrm{diag}}\nolimits}
\def \spf {\mathop{\mathrm{Spf}}\nolimits}
\def \Def {\mathop{\mathrm{Def}}\nolimits}
\def \twist {\mathop{\mathrm{Twist}}\nolimits}
\def \supp {\mathop{\mathrm{Supp}}\nolimits}
\def \Id {{\mathop{\mathrm{Id}}\nolimits}}
\def \bar {\overline}

\def \Ind{\mathop{\mathrm{Ind}}\nolimits}
\def \ind {\mathop{\mathrm{ind}}\nolimits}

\def \mod {\mathop{\mathrm{mod}}\nolimits}
\def \ker {\mathop{\mathrm{Ker}}\nolimits}
\def \coker {\mathop{\mathrm{Coker}}\nolimits}
\def \mult {\mathop{\mathrm{mult}}\nolimits}
\def \vide{\emptyset}
\def \bad {{\mathop{\mathrm{Bad}}\nolimits}}
\def \gal {{\mathop{\mathrm{Gal}}\nolimits}}
\def \Nr {{\mathop{\mathrm{Nr}}\nolimits}}
\def \rn {{\mathop{\mathrm{rn}}\nolimits}}
\def \vol {{\mathop{\mathrm{vol}}\nolimits}}
\def \ad {{\mathop{\mathrm{ad}}\nolimits}}
\def \tr {{\mathop{\mathrm{Tr~}}\nolimits}}
\def \Sp {{\mathop{\mathrm{Sp}}\nolimits}}
\def \lie {{\mathop{\mathrm{Lie}}\nolimits}}
\def \st {{\mathop{\mathrm{St}}\nolimits}}
\def \sp{{\mathop{\mathrm{Sp}}\nolimits}}
\def \card{{\mathop{\mathrm{card}}\nolimits}}
\def \sym{{\mathop{\mathrm{Sym}}\nolimits}}
\def \perv{\mathop{\mathrm{Perv}}\nolimits}
\def \sh {{\mathop{\mathrm{Sh}}\nolimits}}
\def \const {{\mathop{\mathrm{Const}}\nolimits}}
\def \tor {{\mathop{\mathrm{Tor}}\nolimits}}
\def \Zev {{\mathop{\mathrm{Zev}}\nolimits}}

\def \ele{élément }
\def \eles{éléments }
\def \cad{c'est à dire }
\def \rem{{\noindent\textit{Remarque:~}}}
\def \exem{{\noindent \textit{Exemple:~}}}
\def \ssi{~si et seulement si~}
\def \cl {{\mathop{\mathrm{cl}}\nolimits}}
\def \Tw {{\mathop{\mathrm{Tw}}\nolimits}}
\def \ob {{\mathop{\mathrm{Ob}}\nolimits}}
\def \ext {{\mathop{\mathrm{Ext}}\nolimits}}
\def \End {{\mathop{\mathrm{End}}\nolimits}}
\def \inv {{\mathop{\mathrm{inv}}\nolimits}}
\def \fix {{\mathop{\mathrm{Fix}}\nolimits}}

\def\semi{\mathrel{>\!\!\!\triangleleft}}

%
\POP

\setcounter{secnumdepth}{3} \setcounter{tocdepth}{3}

\newcommand{\marque}{\addtocounter{smfthm}{1}
{\smallskip \noindent \textit{\thesmfthm}~---~}}

\renewcommand\atop[2]{\ensuremath{\genfrac..{0pt}{1}{#1}{#2}}}

\title{Réseaux d'induction des représentations elliptiques de Lubin-Tate}

\alttitle{Induced lattices of Lubin-Tate elliptic representations}

\author{Boyer Pascal}
\email{boyer@math.jussieu.fr}
\address{Institut de mathématiques de Jussieu \\ UMR 7586, université Paris 6 \\
175 rue du Chevaleret Paris 13}

\frontmatter

\begin{abstract} Nous étudions la réduction modulo $l$ de certaines représentations elliptiques; pour chacune de ces
représentations nous explicitons un réseau naturellement obtenu par récurrence via l'induction parabolique, en
décrivant le graphe des extensions entre les différents sous-quotients irréductibles de sa réduction modulo $l$. La
motivation essentielle de ce travail est que ces réseaux apparaissent dans la cohomologie des tours de Lubin-Tate.
\end{abstract}

\begin{altabstract} We study the reduction modulo $l$ of some elliptic representations; for each of these
representations, we give a particular lattice naturally obtained by parabolic induction in giving the graph of
extensions between its irreducible sub-quotient of its reduction modulo $l$. The principal motivation for this work, is
that these lattices appear in the cohomology of Lubin-Tate towers.
\end{altabstract}

\subjclass{14G22, 14G35, 11G09, 11G35, 11R39, 14L05, 11G45, 11Fxx}

\keywords{Représentations modulaires, involution de Zelevinski}

\altkeywords{Modular representations, Zelevinski's involution}

\maketitle

\tableofcontents

\pagestyle{headings} \pagenumbering{arabic}

\section*{Introduction}
\renewcommand{\theequation}{\arabic{equation}}
\backmatter

Soit $K$ une extension finie de $\Qm_p$ et $g$ un entier strictement positif. Pour $\pi$ une $\bar
\Qm_l$-représentation irréductible cuspidale entière de $GL_g(K)$, nous étudions dans un premier temps les
sous-quotients irréductibles de la réduction modulo $l$ de la représentation de Steinberg généralisée $\st_s(\pi)$ pour
$s \geq 1$, que l'on note aussi $[\overleftarrow{s-1}]_\pi$ ou $[\overleftarrow{s-1}]$ lorsque $\pi$ est la
représentation triviale de $GL_1(K)$. Pour $s=1$ on sait d'après \cite{vigneras-livre} que $\varrho=r_l(\pi)$ est
irréductible cuspidale mais pas forcément supercuspidale. On note $m(\varrho)$ le cardinal de l'ensemble des $\varrho
\{ i \}$ pour $i \in \Zm$ si celui-ci n'est pas réduit à $\{ \varrho \}$ et sinon $m(\varrho)=l$; le nombre de
sous-quotients irréductibles de $r_l(\st_s(\pi))$ est alors égal au cardinal de
$$\IC_s:=\{ \underline{i}=(i_k)_{k \in \Nm} \in \Nm^\Nm:~s-m(\varrho) \sum_{k=0}^\oo i_k l^k \geq 0 \}.$$
Plus précisément, cf. la proposition \ref{prop-ss-quotient}, à tout $\underline{i} \in \IC_s$ est associé un unique
sous-quotient irréductible de $r_l(\st_s(\pi))$ noté $I_{\underline{i}}([\overleftarrow{s-1}]_\pi)$ caractérisé par le
fait qu'il est de niveau de cuspidalité $\underline{i}=(i_0,\cdots,i_u,0,\cdots)$, i.e., cf. le \S
\ref{superunipotente}, que son image par le foncteur de Jacquet associé au parabolique standard de Levi
$$GL_{s_\varrho(\underline{i})g}(K) \times GL_{m(\varrho)i_kg}(K) \times \cdots \times GL_{m(\varrho)i_ul^u}(K)$$
est cuspidale, où $s_\varrho(\underline{i})=s-m(\varrho)\underline{i}(l)$ avec $\underline{i}(l)=\sum_{k=0}^\oo i_kl^k
\geq 0$. Ces constituants se décrivent en fait simplement à partir de l'involution de Zelevinski $Z_1$, cf. le \S
\ref{invo-Z}, sur les représentations superunipotentes telle qu'elle est définie dans \cite{vigneras-zelevinski} et
dont un calcul explicite est donné dans \cite{ltv}: précisément, proposition \ref{prop-ss-quotient},
$I_{\underline{0}}([\overleftarrow{s-1}])$ est égale à l'image de la représentation triviale de $GL_{sg}(K)$ par $Z_1$ et
$$I_{\underline{i}}([\overleftarrow{s-1}]_\pi) = \Bigl ( I_{\underline{0}}([\overleftarrow{i_0-1}]) \boxtimes \rho_0
\Bigr ) \times \cdots \\ \times \Bigl ( I_{\underline{0}}([\overleftarrow{i_u-1}]) \boxtimes \rho_u \Bigr ) \times
\Bigl ( I_{\underline{0}}([\overleftarrow{s_\varrho(\underline{i})-1}]) \boxtimes \rho_{-1} \Bigr )$$ avec pour
$i=0,\cdots,u$,
$$\rho_i=\st_{m(\varrho)l^i}(\rho\{ \frac{m(\varrho)l^i-s}{2} \} ), \qquad \rho_{-1}=
\varrho \{ \frac{s-s_\varrho(\underline{i})}{2} \},$$ et où pour $\pi_0=<a>$ une représentation superunipotente de
paramètre de Zelevinski $a=(1,r_i)_{1 \leq i \leq k}$, $\pi_0 \boxtimes \varrho$ désigne la représentation de paramètre
de Zelevinski, $(\varrho,r_i)_{1 \leq i \leq k}$, cf. le \S \ref{para-classification}.

Avec les notations du \S \ref{rappel-ql}, en utilisant le morphisme surjectif
$[\overleftarrow{s-2}]_\pi \overrightarrow{\times} [\overleftarrow{0}]_\pi
\twoheadrightarrow [\overleftarrow{s-1}]_\pi$, on construit à la proposition \ref{prop-reseau}, par récurrence, un
réseau $RI_{\bar \Zm_l}([\overleftarrow{s-1}]_\pi)$ de $\st_s(\pi)$ tel que dans le graphe des extensions de sa
réduction modulo $l$, les seules flèches, i.e. les extensions non triviales entre deux constituants irréductibles de
$r_l(\st_s(\pi))$, sont celles qui relient $I_{\underline{i}}([\overleftarrow{s-1}]_\pi)$ et
$I_{\underline{i'}}([\overleftarrow{s-1}]_\pi)$ où $\underline{i} < \underline{i'}$ sont deux éléments successifs de
$\IC_s$ pour l'ordre lexicographique inverse.

On reprend au \S \ref{para-red2} les résultats précédents pour les représentations elliptiques, dites de Lubin-Tate,
$LT_t(s):=[\overleftarrow{t-1},\overrightarrow{s-t}]_\pi$; ainsi, proposition \ref{prop-ss-quotient2}, les constituants
irréductibles de la réduction modulo $l$ de $[\overleftarrow{t-1},\overrightarrow{s-t}]_\pi$ sont indexés par les
éléments $\underline{i} \in \IC_{t-1}$ de telle sorte que le constituant de niveau de cuspidalité $\underline{i}$ est
$$I_{\underline{i}}([\overleftarrow{m(\varrho)\underline{i}(l)-1}]_\pi) \times I_{\underline{0}}(
[\overleftarrow{t_\varrho(\underline{i})-1},\overrightarrow{s-t}]_\pi)$$ où comme précédemment $\underline{i}(l)=
\sum_{k=0}^\oo i_kl^k$ et $t_\varrho(\underline{i})=t-m(\varrho) \underline{i}(l)$. Contrairement au cas $s=t$, l'image
du constituant de niveau de cuspidalité $\underline{i} \in \IC_{t-1}$ de
$r_l([\overleftarrow{t-1},\overrightarrow{s-t}]_\pi)$ par un foncteur de Jacquet n'est, en général, pas irréductible ni
même semi-simple, cf. l'exemple de la fin du paragraphe \ref{para-red2}; ses constituants irréductibles ne sont pas
tous de niveau de cuspidalité $\underline{i}$ à l'exception notable de tout sous-espace et tout quotient irréductible.

En ce qui concerne les réseaux, comme la réduction modulo $l$ de $[\overrightarrow{a}]_\pi$ est irréductible, elle ne
possède à isomorphisme près qu'un seul réseau stable; ainsi les réseaux des Steinberg généralisés construits
précédemment, définissent par induction des réseaux de $[\overleftarrow{t-1},\overrightarrow{s-t}]_\pi$ qui s'avèrent
être tous isomorphes, proposition \ref{prop-compatible} et la remarque qui suit la proposition \ref{prop-reseaubis}; la
description du graphe d'extensions de la réduction modulo $l$ de ces réseaux est similaire à celle de $RI_{\bar
\Zm_l}([\overleftarrow{t-2}]_\pi)$ et donnée à la proposition \ref{prop-reseaubis}.

Dans \cite{boyer-invent2}, nous montrons que les représentations elliptiques du type
$LT_t(\pi,s):=[\overleftarrow{t-1},\overrightarrow{s-t}]_\pi$ apparaissent dans la cohomologie des modèles locaux
définis par Carayol dans \cite{carayol-LT}, généralisant ceux introduits par Lubin et Tate. Dans \cite{boyer-torsion},
nous montrons que les réseaux introduits dans ce papier sont ceux de la $\bar \Zm_l$-cohomologie de ces espaces. En ce
qui concerne les autres représentations elliptiques, la situation semble plus complexe comme le montre l'exemple de la
fin du \S \ref{para-limite}: il n'y a plus, en général, unicité du constituant $\varrho$-superunipotent de la réduction
modulo $l$ d'une représentation elliptique, par contre le cas limite classique suggère que la multiplicité $1$
serait conservée.

\mainmatter

\renewcommand{\theequation}{\arabic{section}.\arabic{subsection}.\arabic{smfthm}}

\section{Rappels sur les représentations de $GL_n(K)$}

Dans la suite $K$ désigne une extension finie de $\Qm_p$ de corps résiduel de cardinal $q$ une puissance de $p$. Dans
cet article nous étudierons des représentations admissibles $\pi$ de $GL_d(K)$ à coefficients dans $\bar \Qm_l$, $\bar
\Zm_l$ ou $\bar \Fm_l$ où $l$ est un nombre premier distinct de $p$. En ce qui concerne la torsion, suivant
\cite{boyer-invent2}, nous la noterons $\pi\{ n \}$ de sorte que l'action d'un élément $g \in GL_d(K)$ est donnée par
$\pi(g) |\det g|^n$ où $|-|$ est la valeur absolue sur $K$.

\subsection{Induction et foncteur de Jacquet}

\begin{defi} \label{defi-parab}
Pour une suite $r_1, r_2, \cdots ,r_k$ telle que $\sum_{i=1}^k r_i=d$, on note $P_{r_1,r_2,\cdots,r_k}$ le sous-groupe
parabolique de $GL_d$ standard associé au sous-groupe de Levi $GL_{r_1}(F_v) \times GL_{r_2}(F_v) \times \cdots \times
GL_{r_k}(F_v)$ et on note $N_{r_1,\cdots,r_k}$ son radical unipotent. On mettra en exposant $op$ pour désigner les
paraboliques opposés. Pour $\underline{\lambda}=(d=\lambda_1 \geq \lambda_2 \geq \cdots \geq \lambda_r)$, on note
$P_{\underline{\lambda}}=P_{\lambda_r,\lambda_{r-1}-\lambda_r,\cdots,\lambda_1-\lambda_2}$.
\begin{itemize}
\item Pour $\pi_1$ et $\pi_2$ des représentations de respectivement $GL_{n_1}(K)$ et $GL_{n_2}(K)$, $\pi_1 \times \pi_2$
désigne l'induite à $GL_{n_1+n_2}(K)$
$$\Ind_{P_{n_1,n_1+n_2}(K)} \pi_1\{ n_2/2 \} \otimes \pi_2 \{ -n_1/2\}$$
relativement au parabolique standard, que l'on notera aussi parfois sous sa forme normalisée $\ind_{P_{n_1,n_2}(K)}
\pi_1 \otimes \pi_2$. Dans le cas du parabolique opposé au standard, on notera l'induite correspondante
$\pi_1 \times_{op} \pi_2$.

\item \textit{Foncteurs de Jacquet}: soit $P=MN$ un parabolique de $GL_d$ de Lévi $M$ et de radical unipotent $N$. Pour
$\pi$ une représentation admissible de $GL_d(K)$, l'espace des vecteurs $N(K)$-coinvariants est stable sous l'action de
$M(K) \simeq P(K)/N(K)$. On notera $J_P(\pi)$ cette représentation tordue par $\delta_P^{-1/2}$: c'est un foncteur
exact.
\end{itemize}
\end{defi}

\textit{Formules d'adjonction}: $J_P$ est un adjoint à gauche de $\ind_P$, c'est \textit{la réciprocité de Frobenius}.
Par ailleurs \textit{la deuxième formule d'adjonction}, prouvée par Dat dans \cite{dat-finitude} en toute généralité,
affirme que $\times_{P^{op}} \otimes \delta_P^{-1}$ est un adjoint à gauche de $J_P$.

\subsection{Représentations elliptiques sur $\bar \Qm_l$}
\label{rappel-ql}

On rappelle dans ce paragraphe, quelques notations sur les $\bar \Qm_l$-représentations admissibles de
$GL_d(K)$, tirées de \cite{boyer-invent2}.

\begin{defis} Soit $g$ un diviseur de $d=sg$ et $\pi$ une représentation cuspidale irréductible de
$GL_g(K)$:

- les sous-quotients irréductibles de $V(\pi,s):=\pi\{ \frac{1-s}{2}\} \times \pi\{ \frac{3-s}{2} \} \times
\cdots \times \pi\{ \frac{s-1}{2}\}$ seront dits elliptiques de type $\pi$;

- $V(\pi,s)$ possède un unique quotient (resp. sous-espace) irréductible que l'on notera
$[\overleftarrow{s-1}]_{\pi}$ (resp.  $[\overrightarrow{s-1}]_{\pi})$); c'est une représentation de
Steinberg (resp. de Speh) généralisée notée habituellement $\st_s(\pi)$ (resp. $\speh_s(\pi)$);

- pour $\pi_1$ et $\pi_2$ des représentations respectivement de $GL_{t_1g}(K)$ et $GL_{t_2g}(K)$, on
notera $\pi_1 \overrightarrow{\times} \pi_2$ (resp. $\pi_1 \overleftarrow{\times} \pi_2$) l'induite
$\pi_1\{ -t_2/2 \} \times \pi_2\{ t_1/2 \}$ (resp. $\pi_1\{ t_2/2 \} \times \pi_2\{ -t_1/2 \}$)
relativement au parabolique
standard, où on omet la référence à $g$ car en général le contexte sera clair.
On notera comme précédemment $\overrightarrow{\times}_{op}$ et $\overleftarrow{\times}_{op}$ les
mêmes induites relativement au parabolique opposé.
\end{defis}

\noindent \textbf{On conviendra} que $\Pi_t \times [\overleftrightarrow{-1}]_{\pi}:= \Pi_t$.

\medskip

\noindent \textit{Propriétés:}
\begin{itemize}
\item[-] pour $(\pi_i)_{1 \leq i \leq 3}$ des représentations de $GL_{t_ig}(K)$, on a les égalités suivantes:
$$\begin{array}{rl}
(\pi_1 \overrightarrow{\times} \pi_2)^\vee \simeq & \pi_1^\vee \overleftarrow{\times} \pi_2^\vee \\
(\pi_1 \overrightarrow{\times} \pi_2) \overrightarrow{\times} \pi_3 = & \pi_1 \overrightarrow{\times}
(\pi_2
\overrightarrow{\times} \pi_3) \\
(\pi_1 \overleftarrow{\times} \pi_2) \overleftarrow{\times} \pi_3= & \pi_1 \overleftarrow{\times} (\pi_2
\overleftarrow{\times} \pi_3)
\end{array}$$
On a les mêmes égalités pour les versions relativement aux paraboliques opposés au standard. En outre si
$\pi_1$ et $\pi_2$ sont elliptiques de type $\pi$, il en est de même de $\pi_1 \overrightarrow{\times}
\pi_2$ et donc de $\pi_1 \overleftarrow{\times} \pi_2$.

\item[-] soit $g$ un diviseur de $d=sg$ et $\pi$ une représentation irréductible cuspidale de $GL_g(K)$, on a
$$J_{P_{tg,(s-t)g}}([\overleftarrow{s-1}]_{\pi})=[\overleftarrow{t-1}]_{\pi\{ (s-t)/2\} } \otimes
[\overleftarrow{s-t-1}]_{\pi\{ -t/2 \} }$$
$$J_{P_{tg,(s-t)g}}([\overrightarrow{s-1}]_{\pi})=[\overrightarrow{t-1}]_{\pi\{ (t-s)/2\} } \otimes
[\overrightarrow{s-t-1}]_{\pi\{ t/2\} }$$
$$J_{P_{tg,(s-t)g}^{op}}([\overleftarrow{s-1}]_{\pi})=[\overleftarrow{t-1}]_{\pi\{ (t-s)/2\} } \otimes
[\overleftarrow{s-t-1}]_{\pi\{ t/2 \}}$$
$$J_{P_{tg,(s-t)g}^{op}}([\overrightarrow{s-1}]_{\pi})=[\overrightarrow{t-1}]_{\pi\{ (s-t)/2\} } \otimes
[\overrightarrow{s-t-1}]_{\pi\{ -t/2\} }$$

\item[-] l'induite
$$[\overleftarrow{t-1}]_{\pi} \overrightarrow{\times} [\overrightarrow{s-t-1}]_{\pi} \qquad (\hbox{resp. }
[\overleftarrow{t-1}]_{\pi} \overrightarrow{\times}_{op} [\overrightarrow{s-t-1}]_{\pi})$$ est de
longueur $2$; on notera
$[\overleftarrow{t-1},\overrightarrow{s-t}]_{\pi}$ (resp. $[\overrightarrow{s-t},\overleftarrow{t-1}]_{\pi})$)
son unique sous-espace irréductible, et
$[\overleftarrow{t},\overrightarrow{s-t-1}]_{\pi}$ (resp. $[\overrightarrow{s-t-1},\overleftarrow{t}]_{\pi})$)
son unique quotient irréductible;

\item[-] par dualité, l'induite
$$[\overleftarrow{t-1}]_{\pi} \overleftarrow{\times} [\overrightarrow{s-t-1}]_{\pi} \qquad (\hbox{resp. }
[\overleftarrow{t-1}]_{\pi} \overleftarrow{\times}_{op} [\overrightarrow{s-t-1}]_{\pi})$$ est de longueur
$2$ avec
$[\overrightarrow{s-t-1},\overleftarrow{t}]_{\pi}$ (resp. $[\overleftarrow{t},\overrightarrow{s-t-1}]_{\pi})$)
pour unique sous-espace irréductible et
$[\overrightarrow{s-t},\overleftarrow{t-1}]_{\pi}$ (resp. $[\overleftarrow{t-1},\overrightarrow{s-t}]_{\pi})$)
pour unique quotient irréductible.
\end{itemize}

\rem dans \cite{boyer-invent2}, nous montrons que les représentations $[\overleftarrow{t-1},\overrightarrow{s-t}]_\pi$
apparaissent dans la cohomologie des tours dites de Lubin-Tate définis dans \cite{carayol-LT}; pour cela on les notera
sous la forme $LT_t(\pi,s)$.

\subsection{Classification à la Zelevinski sur $\bar \Fm_l$}
\label{para-classification}

On rappelle que $l$ et $p$ désignent des nombres premiers distincts et que $q$ est une puissance de $p$. On note
$e_l(q)$ l'ordre de l'image de $q$ dans $\Fm_l^\times$. On dit que $l$ est \textit{banal} pour $GL_d(K)$ si $e_l(q)>d$.

\begin{defi} Une représentation $\varrho$ de $GL_d(K)$ est dite cuspidale si pour tout parabolique
propre $P$, $J_{P}(\pi)$ est nul. Elle sera dite supercuspidale si elle n'est pas un sous-quotient d'une induite
parabolique propre. On notera $\cusp_l=\bigcup_{n \geq 1} \cusp_l(n)$ la réunion de l'ensemble des classes
d'isomorphismes des $\bar \Fm_l$-représentations irréductibles cuspidales de $GL_n(K)$.
\end{defi}

\marque \textit{Segments de Zelevinski}: le caractère non ramifié $\nu: g \in GL_n(K) \mapsto p^{\val (\det g)} \in
\bar \Fm_l$ agit sur $\cusp_l$ par multiplication.

\begin{defis} \label{defi-parametre}
Soit $\varrho$ une $\bar \Fm_l$-représentation irréductible cuspidale de $GL_g(K)$:
\begin{itemize}
\item la droite associée à $\varrho$ est par définition l'ensemble
$\{ \varrho\{ i\} ~/~ i \in \Zm \}$; il est de cardinal fini $\e(\varrho)$ un diviseur de $e_l(q)$,
cf. \cite{vigneras-induced} p.51.
On pose comme dans loc. cit. $m(\varrho)=\e(\varrho)$ si $\e(\varrho)>1$ et sinon $m(\varrho)=l$;

\item un segment de Zelevinski associé à $\varrho$ et de longueur $r \geq 1$ est une suite $$(\varrho,r)=(\varrho,\varrho\{ 1 \},\cdots,\varrho\{r-1 \} ).$$

\item Un \textbf{paramètre de Zelevinski} est un multi-ensemble $a=(\varrho_i,r_i)_{1 \leq i \leq t}$
de segments de Zelevinski dans $\cusp_l$; le support $s$ de $a$ est la réunion des multi-ensembles
$(\varrho_i,r_i)$.

\item Un paramètre de Zelevinski de la forme $(\varrho_i,r)_{1 \leq i \leq t}$ où
$(\varrho_1,\cdots,\varrho_t)$
est un segment $(\varrho,m(\varrho)l^u)$ pour $u \geq 0$, est appelé \textbf{un cycle}.

\item Un paramètre de Zelevinski qui ne contient pas de cycle est dit \textbf{restreint}. Si son support
cuspidal est supercuspidal le paramètre est dit \textbf{supercuspidal}.
\end{itemize}
\end{defis}

\rem (cf. \cite{vigneras-induced} p.54) étant donné un paramètre de Zelevinski $a$, si le segment
$(\varrho,r)$ apparaît dans $a$ avec $a$ cuspidal mais non supercuspidal, alors il existe une
représentation $\varrho'$ supercuspidale et un entier $s=m(\varrho') l^u$ tel que
$\varrho=\st_s(\varrho')$. En remplaçant
systématiquement de tels $(\varrho,r)$ par le cycle $(\varrho' \nu^i,r)_{0 \leq i < s}$, on obtient un
paramètre supercuspidal $a_{sc}$. En procédant en sens inverse on peut remplacer tous les cycles
par des segments de la forme $(\st_s(\varrho)\{ i \} ,r)_{0 \leq i <s}$ afin d'obtenir un paramètre
restreint $a_{rst}$.

\marque \textit{Classification des représentations irréductibles,} cf. \cite{vigneras-induced} V.9:
étant donné un paramètre de Zelevinski $a=(\varrho_i,r_i)_{1 \leq i \leq t}$, on lui associe une
représentation irréductible $<a>$ qui vérifie les propriétés suivantes:
\begin{itemize}
\item si $a=(\varrho,r)$ alors $<(\varrho,r)>$ est l'unique, à isomorphisme près\footnote{i.e. la
multiplicité n'est pas forcément égale à $1$ dans le cas non banal}, sous-représentation de
l'induite parabolique $\varrho \times \varrho \nu \times \cdots \times \varrho \nu^{r-1}$ telle que
$$J_{P_{n,n,\cdots,n}}(<(\varrho,r)>)=\varrho \times \cdots \otimes \varrho \{ r-1\};$$

\item de manière générale $<a>$ est l'unique sous-quotient de
$$\pi(a):=<(\varrho_1,r_1)> \times \cdots \times <(\varrho_t,r_t)>$$
dont le niveau de Whittaker, cf. loc. cit., est le même que celui de $\pi(a)$;

\item on a $<a>=<a_{sc}>=<a_{rst}>$ et le support supercuspidal (resp. cuspidal) de $<a>$ est le
support de $a_{sc}$ (resp. de $a_{rst}$);

\item pour $a$ et $a'$ deux paramètres de Zelevinski supercuspidaux distincts alors $<a>$ et $<a'>$
ne sont pas isomorphes;

\item l'application $a \mapsto <a>$ est surjective.
\end{itemize}

\begin{prop} (\cite{vigneras-induced} III.5.14) L'induite parabolique
$$\varrho \overrightarrow{\times} \cdots \overrightarrow{\times} \varrho = \varrho\{\frac{1-s}{2} \} \times
\cdots \times \varrho \{ \frac{s-1}{2} \}$$ admet un unique sous-quotient non dégénéré que l'on note $\st_s(\varrho)$.
La représentation $\st_s(\varrho)$ est cuspidale si et seulement si
$$s=1,m(\varrho),m(\varrho)l,\cdots,m(\varrho) l^u,\cdots$$
La réunion de ces dernières avec les supercuspidales forment l'ensemble des représentations
cuspidales.
\end{prop}

\subsection{Réduction modulo $l$ des $\bar \Qm_l$-représentations entières: généralités}
\label{red-elliptique}

Une $\bar \Qm_l$-représentation lisse de longueur finie $\pi$ de $GL_n(K)$ est dite entière s'il
existe une extension finie $E/ \Qm_l$ d'anneau des entiers $\OC_E$ et une $\OC_E$-représentation
$L$ de $GL_n(K)$ qui est un $\OC_E$-module libre tel que $\bar \Qm_l \otimes_{\OC_E} L \simeq \pi$
et tel que $L$ est un $\OC_E GL_n(K)$-module de type fini. Soit $\k_E$ le corps résiduel de
$\OC_E$, on dit que $\k_E \otimes_{\OC_E} L$ est la réduction de L et que $\bar \Fm_l
\otimes_{\OC_E} L$ est la réduction modulo $l$ de $L$.

\marque \textit{Principe de Brauer-Nesbitt}: la semi-simplifiée de $\bar \Fm_l \otimes_{\OC_E} L$ est une $\bar
\Fm_l$-représentation de $GL_n(K)$ de longueur finie qui ne dépend pas du choix de $L$. Son image dans le groupe de
Grothendieck sera notée $r_l(\pi)$ et dite la réduction modulo $l$ de $\pi$.

\rem la réduction modulo $l$, commute avec l'induction parabolique et les foncteurs de Jacquet. Elle ne respecte pas le
caractère supercuspidal mais la réduction modulo $l$ d'une représentation irréductible entière supercuspidale est
irréductible cuspidale. Ainsi, pour $\pi$ une $\bar \Qm_l$-représentation irréductible cuspidale entière, on dira
qu'une $\bar \Fm_l$-représentation est $\pi$-unipotente si elle est $r_l(\pi)$-superunipotente.

On rappelle, cf. \cite{vigneras-livre} II.4.12, qu'une $\bar \Qm_l$-représentation irréductible
supercuspidale de $GL_n(K)$ est entière si et seulement si son caractère central est entier. De
manière générale une représentation irréductible de $GL_n(K)$ est entière si et seulement si les
constituants irréductibles de son support supercuspidal $sc(\pi)$ sont entiers.

\marque Étant donnée une représentation $\Pi$ de $GL_d(K)$, on lui associe le graphe orienté $\Gamma_\Pi$ défini comme
suit:
\begin{itemize}
\item ses sommets sont les différents facteurs irréductibles, $\pi_1,\cdots,\pi_c$, de $\Pi$ vu dans le groupe de
Grothendieck;

\item une flèche relie $\pi_i$ à $\pi_j$ si et seulement s'il existe un sous-quotient de $\Pi$
qui est une extension \textit{non triviale} de $\bar \pi_j$ par $\bar \pi_i$.
\end{itemize}

Dans le cas où $\Pi$ est une $\bar \Qm_l$-représentation irréductible entière de $GL_d(K)$, pour un réseau stable
$\Lambda$, $\Gamma_\Lambda$ désigne, à la manière de \cite{bellaiche-ribet}, le graphe de la représentation $\Lambda
\otimes_{\bar \Zm_l} \bar \Fm_l$ au sens ci-dessus. Dans les paragraphes suivants nous construirons des exemples de
réseaux pour les représentations elliptiques $[\overleftarrow{t-1},\overrightarrow{s-t}]_\pi$ où $\pi$ est une $\bar
\Qm_l$-représentation irréductible cuspidale; ce sont ceux qui apparaissent \og{} naturellement \fg{} dans la
cohomologie des tours de Lubin-Tate, cf. \cite{boyer-torsion}. Signalons par ailleurs que l'on dispose toujours du
réseau particulier donné par le lemme suivant et dont le graphe ne possède aucune flèche.

\nocite{bellaiche-ribet2}

\begin{prop} \label{prop-ss}
(cf. le lemme 6.11 de \cite{dat-duke}) Étant donnée une $\bar \Qm_l$-représentation entière
irréductible $\pi$ de $GL_n(K)$, il existe un réseau $\pi^l$ tel que $\pi^l \otimes_{\bar \Zm_l}
\bar \Fm_l$ est semi-simple.
\end{prop}

\section{Représentations $\varrho$-superunipotentes}
\label{superunipotente}

\subsection{Définitions}

Rappelons qu'un bloc dans $\mod_{\bar \Fm_l} GL_n$ est une sous-catégorie abélienne pleine qui est
un facteur direct sans être le produit direct de deux sous-catégories abéliennes non nulles.
Soit $(M,\rho)$ une $\bar \Fm_l$-représentation irréductible cuspidale d'un sous-groupe de Levi
$M$ d'un parabolique de $GL_n$. On considère
$$\irr_{\bar \Fm_l GL_n,[\rho]}, \qquad \BF_{\bar \Fm_l GL_n,[\rho]}$$
avec $\irr_{\bar \Fm_l GL_n,[\rho]}$ l'ensemble des classes d'isomorphismes des $\bar
\Fm_l$-représentations irréductibles de $GL_n$ qui sont des sous-quotients de $\Ind_{M'}^{GL_n}
\rho'$ pour $(M',\rho')$ inertiellement équivalent à $(M,\rho)$ et $\BF_{\bar \Fm_l GL_n,[\rho]}$
la sous-catégorie abélienne de $\mod_{\bar \Fm_l} GL_n$ constituée des $\bar
\Fm_l$-représentations de $GL_n$ dont tous les sous-quotients irréductibles sont dans $\irr_{\bar
\Fm_l GL_n,[\rho]}$.

\begin{defi} Dans le cas où $\rho$ est la représentation triviale du groupe des matrices diagonales, le
bloc correspondant est appelé le bloc unipotent. Les représentations irréductibles de ce bloc sont appelés les
\textit{unipotentes}. Les représentations \textit{superunipotentes} sont les unipotentes caractérisées par l'une des
deux propriétés équivalentes suivantes:
\begin{itemize}
\item il existe un vecteur fixe par le sous-groupe d'Iwahori standard;

\item la restriction parabolique au sous-groupe des matrices diagonales, n'est pas nulle.
\end{itemize}
\end{defi}

\begin{theo} (cf. \cite{vigneras-induced} théorème III.6)
La catégorie $\mod_{\bar \Fm_l} GL_n$ est un produit direct de blocs, chaque bloc étant de la
forme $\BF_{\bar \Fm_l GL_n,[\rho]}$ pour $\rho$ irréductible supercuspidal.
\end{theo}

\begin{prop} \label{prop-superunip}
(cf. \cite{vigneras-induced} IV.2.5, IV.6.2, IV.6.3) Pour $\rho$ irréductible cuspidal, il existe un
groupe $G'$ de la forme $\prod_i GL_{n_i}(E_i)$ avec $E_i$ une extension finie de $K$, ainsi
qu'une bijection entre l'ensemble des représentations irréductibles cuspidal unipotentes de $G'$
et $\irr_{\bar \Fm_l GL_n,[\rho]}$ qui respecte le support cuspidal.
\end{prop}

\rem dans le cas banal, on a même une équivalence de catégorie entre le bloc unipotent de $G'$ et
$\BF_{\bar \Fm_l GL_n,[\rho]}$.

\begin{defi} \label{defi-rho}
Soit $\varrho$ une représentation cuspidale de $GL_g(K)$ pour $g$ un diviseur de $d=sg$, et soit
$$\rho=\varrho\{\frac{(s-1)(g-1)}{2}\} \otimes \cdots \otimes \varrho\{ \frac{(1-s)(g-1)}{2}\}$$
la représentation cuspidale de $M=GL_g(K)^s$. On dira d'une représentation dans le bloc $\BF_{\bar \Fm_l GL_n,[\rho]}$,
qu'elle est $\varrho$-superunipotente si sa restriction parabolique au groupe diagonal par blocs $M$ est non nulle.
\end{defi}

\rem les représentations superunipotentes sont avec cette définition les représentations
$1_K$-superunipotentes, où $1_K$ désigne la représentation triviale de $K^\times$.

\rem les foncteurs de Jacquet et d'induction ne respectent pas le caractère $\varrho$-superunipotent.

\begin{defi} (\cite{vigneras-zelevinski} \S 4.4)
Pour $\varrho$ une représentation cuspidale de $GL_g(K)$ et $g$ un diviseur de $d=sg$, le groupe de Grothendieck
$\varrho$-superunipotent $\groth_{\varrho-su}$ est le quotient du groupe de Grothendieck $\groth_{\varrho-u}$ des
représentations $\varrho$-unipotentes par les non $\varrho$-superunipotentes. On désigne alors par
$\MC_{\varrho}:\groth_{\varrho-u} \twoheadrightarrow \groth_{\varrho-su}$ la projection  associée.
\end{defi}

\begin{coro} \label{coro-superunipotent}
(cf. \cite{vigneras-induced} IV.6.3 p47) Pour $\varrho$ irréductible cuspidale et
$\rho$ comme dans la définition \ref{defi-rho},
on a une bijection $F_{G,G'}$ entre les
représentations $\varrho$-superunipotentes du bloc $\BF_{\bar \Fm_l GL_n,[\rho]}$ et les
représentations superunipotentes de $GL_s(E)$ pour une certaine extension finie $E/K$. Ces
$F_{G,G'}$ commutent alors aux foncteurs de Jacquet et d'induction au sens où si $P$ est un parabolique
de $G$ de Lévi $M$ associé à un parabolique $P'$ de $G'$ de Lévi $M'$ alors
$$F_{G,G'} \circ \Ind_{P'}^{G'}=\Ind_P^G \circ F_{G,G'}, \qquad F_{M,M'} \circ J_{M'}^{G'}=J_P^G \circ
F_{G,G'}.$$
\end{coro}

\begin{proof} Le fait qu'il n'y ait qu'un seul $GL_s(E)$ découle de la construction de loc. cit. IV.3.3.
Le
groupe de Grothendieck superunipotent correspond à celui des modules sur l'algèbre de
Hecke-Iwahori. La bijection découle alors de l'isomorphisme entre ces deux algèbres et la
commutativité des foncteurs de Jacquet et d'induction, modulo les non superunipotentes, découle de
loc. cit. p.47.
\end{proof}

\begin{defi} \label{defi-boxtimes}
Avec les notations du corollaire précédent, pour $\pi_0$ une représentation irréductible
superunipotente de $G'=GL_s(E)$, on note $\pi_0 \boxtimes \varrho$ la représentation
$\varrho$-superunipotente $F_{G,G'}(\pi_0)$.
\end{defi}

\rem si $\pi_0=<a>$ pour un paramètre de Zelevinski $a=(1,r_i)_{1 \leq i \leq k}$ alors
$\pi_0 \boxtimes \varrho$ a pour paramètre de Zelevinski, $(\varrho,r_i)_{1 \leq i \leq k}$.

\subsection{Niveau de cuspidalité}

\begin{prop} \label{niveau-cusp}
Étant donnée une représentation cuspidale $\varrho$ de $GL_g(K)$, pour toute représentation $\pi$ de $GL_{sg}(K)$, du
bloc $\varrho$-unipotent, on considère les paraboliques $P$ de $GL_d$ tels que le foncteur de Jacquet $J_P$ appliqué à
$\pi$ soit non nul. Les éléments minimaux pour l'inclusion de cet ensemble de parabolique ont alors des sous-groupes de
Lévi conjugués par des matrices de permutation. L'ordre de Bruhat sur les partitions de $s$ donne alors une relation
d'ordre partiel sur le bloc $\rho$-unipotent de $GL_{sg}(K)$: on y référera comme \textbf{le niveau de cuspidalité}
relativement à $\varrho$.
\end{prop}

\begin{proof}
Soit $P$ un élément minimal de Lévi $\prod_{i=1}^r GL_{n_i}(K)$; d'après la transitivité des foncteurs de Jacquet on en
déduit que $J_P(\pi)$ est cuspidal. D'après la réciprocité de Frobenius $\pi$ est un sous-espace d'une induite
parabolique $\pi_1 \times \cdots \times \pi_r$, où les $\pi_i$ sont des représentations irréductibles cuspidales. Le
résultat découle alors de l'exactitude à gauche du foncteur de Jacquet et du fait que la proposition est classiquement
vérifiée pour $\pi_1 \times \cdots \times \pi_r$ où les $\pi_i$ sont cuspidaux.
\end{proof}

\rem d'après l'unicité du support supercuspidal, la partition de $s$ associée au niveau de cuspidalité d'une
représentation $\rho$-unipotente est de la forme $\underline{\lambda}=(\lambda_1 \geq \cdots \lambda_r)$ où pour tout
$1 \leq i < r$, $\lambda_{i}-\lambda_{i+1}$ est soit égal à $1$ ou de la forme $m(\varrho)l^k$; pour tout $k \geq 0$
notons $m_k$ (resp. $m_{-1}$) le nombre d'indice $1 \leq i <r$ tels que $\lambda_i-\lambda_{i+1}=m(\varrho)l^k$ (resp.
$\lambda_i-\lambda_{i+1}=1$) avec donc
$$s=m_{-1}+m_0 m(\varrho)+m_1 m(\varrho)l+ \cdots +m_u m(\varrho)l^u, \qquad 0 \leq m_{-1}, \cdots,m_u.$$
Les $\varrho$-superunipotentes correspondent à $m_{-1}=s$ et $m_k=0$ pour tout $k \geq 0$. On notera \textit{le niveau
de cuspidalité} sous la forme $(m_0,\cdots,m_u)$ et \textit{l'ordre de Bruhat} est alors l'ordre lexicographique
inverse. Par ailleurs, on notera que pour $\pi$ irréductible le niveau cuspidal de tout constituant irréductible de
$J_P(\pi)$ est supérieur ou égal à celui de $\pi$; le même phénomène se produit pour l'induction parabolique.

\noindent \textit{Exemples}: on fixe une représentation $\varrho$ irréductible cuspidale de $GL_g(K)$ et pour
$i=0,\cdots,u$, on pose $\rho_i=\st_{m(\varrho)l^i} (\varrho\{ \frac{m(\varrho)l^i-s}{2}\})$ et $\rho_{-1}=\varrho$: ce
sont des représentations cuspidales.
\begin{itemize}
\item Soit $s=m_{-1}+m_0 m(\varrho)+m_1 m(\varrho)l+ \cdots +m_u m(\varrho)l^u$ avec
$0 \leq m_{-1} < m(\varrho)$ et $0 \leq m_0,\cdots,m_u < l$; le niveau de cuspidalité de $\st_s(\varrho)$ est alors
$(m_0,\cdots,m_u)$ avec
$$\st_s(\varrho) \simeq \st_{m_{-1}}(\varrho_{-1}) \times \st_{m_0} (\rho_{0}) \times \cdots \times \st_{m_u} (\rho_u).$$

\item Pour $\pi_1,\cdots,\pi_r$
des représentation irréductibles $\varrho$-unipotentes de respectivement $GL_{n_ig}$ et de niveau de cuspidalité
$(m_{0}(i),\cdots,m_u(i))$; le niveau de cuspidalité de $\pi_1 \otimes \cdots \otimes \pi_r$ est alors égal à
$\sum_{i=1}^r \Bigl ( m_{0}(i),\cdots,m_u(i) \Bigr ) $.

\item Soient pour $i=-1,0,\cdots,u$ des entiers $m_i \geq 0$ et des représentations $\pi_i$
irréductibles superunipotentes alors
$$\Bigl ( \pi_{-1} \boxtimes \rho_{-1} \Bigr ) \times \cdots \times \Bigl ( \pi_u \boxtimes \rho_u
\Bigr )$$ est irréductible, cf. \cite{vigneras-induced} \S V.3, de niveau de cuspidalité $(m_{0},\cdots,m_u)$.
\end{itemize}

\subsection{Involution de Zelevinski}
\label{invo-Z}

Dans ce paragraphe $\varrho$ désigne une représentation irréductible cuspidale de $GL_g(K)$; nous allons reprendre les
résultats principaux de \cite{vigneras-zelevinski} pour le bloc $\varrho$-unipotent. Pour tout parabolique standard $Q$
de $GL_d$, on note, pour $i \in \Nm$, $\PC_i(Q)$ l'ensemble des paraboliques standard de $Q$ dont les Levi ont $i+1$
facteurs $GL_{r_k}$ de plus que $Q$; pour $Q=GL_d$, on notera simplement $\PC_i$. Pour toute représentation $\rho$ du
bloc $\varrho$-unipotent, on a alors un complexe augmenté, \cite{vigneras-zelevinski} \S 3.8:
$$0 \rightarrow \rho \rightarrow \bigoplus_{P \in \PC_0} T_P (\rho) \rightarrow \bigoplus_{P \in \PC_1} T_P(\rho)
\rightarrow \cdots \rightarrow \bigoplus_{P \in \PC_{s-1}} T_P (\rho)$$
où $T_P=\ind_P^G \circ J_P$, dont l'homologie $H^i(K^\bullet(\rho))$ vérifie les points suivants,
cf. \cite{vigneras-zelevinski} théorème 4.6:
\begin{itemize}
\item[(i)] $H^{s-1}(K^\bullet(\rho))$ est non nul si et seulement si $\rho$ est $\varrho$-superunipotente auquel cas,
pour $\rho$ irréductible, il admet un unique quotient irréductible qui est $\varrho$-superunipotent et que
l'on note $Z_\varrho(\rho)$;

\item[(ii)] pour $\rho$ une représentation irréductible
$\varrho$-superunipotente, $Z_\varrho(\rho)$ est le seul constituant $\varrho$-superunipotent de $\bigoplus_i
H^i(K^\bullet(\rho))$. L'application $\rho \mapsto Z_\varrho(\rho)$ est un involution qui préserve l'irréductibilité.
\end{itemize}

\begin{lemm} \label{lem-Zcommute}
L'involution $Z_{\varrho}$ défini dans le groupe de Grothendieck des représentations $\varrho$-superunipotente
commute aux foncteurs de Jacquet et à l'induction.
\end{lemm}

\begin{proof} Rappelons que, cf. \cite{ze} corollaire du \S 1.3,
pour deux partitions $\underline{\beta}=(\beta_1 \leq \cdots \leq \beta_r=d)$ et
$\underline{\gamma}=(\gamma_1 \leq  \cdots \leq \gamma_s=d)$ et pour toute représentation $\pi$ de
$GL_{\beta_1}(K) \times GL_{\beta_2-\beta_1}(K) \times \cdots \times GL_{\beta_r-\beta_{r-1}}(K)$,
la représentation $J_{P_{\underline{\gamma}}}^{GL_d} \circ \ind_{P_{\underline{\beta}}}^{GL_d}(\pi)$ admet
$\ind_{P_{\underline{\beta} \cap \underline{\gamma}}}^{P_{\underline{\gamma}}} \circ J_{P_{\underline{\beta}
\cap \underline{\gamma}}}^{P_{\underline{\beta}}}(\pi)$ comme quotient.
Ainsi pour tout parabolique $Q$, on en déduit avec les notations ci-dessus, des surjections
pour tout $i=0,\cdots,s-1$
$$J_Q \Bigl ( \bigoplus_{P \in \PC_i} T_P(\pi) \Bigr ) \twoheadrightarrow \bigoplus_{P \in \PC_i(Q)}
T_P \Bigl ( J_Q^{GL_d} (\pi) \Bigr ), \qquad \bigoplus_{P \in \PC_i} T_P \Bigl ( \ind_Q^{GL_d} \pi \Bigr )
\twoheadrightarrow \ind_{Q}^{GL_d} \Bigl ( \bigoplus_{P \in \PC_i(Q)} T_P(\pi) \Bigr )$$
Par ailleurs pour $\rho$ une représentation $\varrho$-superunipotente, on a $T_P(\rho)=0$ pour tout
$P \in \PC_i$ si $i > s-1$, de sorte que l'on obtient des surjections
$$J_Q^{GL_d} H^{s-1}(K^\bullet(\pi)) \twoheadrightarrow H^{s-1}(K^\bullet(J_Q^{GL_d}(\pi))),
\qquad H^{s-1}(K^\bullet(\ind_Q^{GL_d} \pi)) \twoheadrightarrow \ind_Q^{GL_d} H^{s-1}(K^\bullet(\pi)).$$
Le résultat découle alors des faits suivants:
\begin{itemize}
\item d'après (ii) ci-dessus, les parties $\varrho$-superunipotentes de
$H^i(K^\bullet(J_Q^{GL_d}(\pi)))$, $J_Q^{GL_d} H^{s-1}(K^\bullet(\pi))$ et $H^i(K^\bullet(\ind_Q^{GL_d}(\pi)))$
sont nuls pour $i \neq s-1$;

\item on a égalité des sommes alternées
$$J_Q^{GL_d} H^{*}(K^\bullet(\pi)) = H^{*}(K^\bullet(J_Q^{GL_d}(\pi))),
\qquad H^{*}(K^\bullet(\ind_Q^{GL_d} \pi)) = \ind_Q^{GL_d} H^{*}(K^\bullet(\pi)).$$
\end{itemize}
\end{proof}

Dans \cite{ltv} les auteurs décrivent l'involution de Zelevinski $Z_{\bar \Fm_l}$ sur les
représentations \textbf{superunipotentes} ce qui d'après \ref{coro-superunipotent} fourni le
calcul de $Z_{\varrho}$ sur les représentation $\varrho$-superunipotentes pour toute représentation
$\varrho$ irréductible cuspidale. En particulier le résultat suivant découle de leur description.

\begin{lemm} \label{lem-Ztrivial}
Pour $\epsilon(\rho)>2$ (resp. $\epsilon(\rho)=2$), $Z_{\rho}(<(\rho,s)>)$ est égale à $<a>$
pour le paramètre de Zelevinski $a$:
$$(\rho \nu^{-q},q+1),~(\rho \nu^{1-q},q+1),\cdots,(\rho \nu^{r-1-q},q+1),~(\rho \nu^{r+1-q},q),~(\rho
\nu^{r+2-q},q),\cdots,(\rho \nu^{-1-q},q)$$
où $r$ est le reste de la division euclidienne de $s$ par $\epsilon(\rho)-1>1$:
$s=q(\epsilon(\rho)-1)+r$, (resp. $<(\rho,s)>$ pour $s$ pair et $<\rho(1),s>$ pour $s$ impair).
\end{lemm}

\subsection{Cas limite classique}
\label{para-limite}

Cela correspond à $q \equiv 1 \mod l$ et $d<l$. On s'intéresse aux représentations superunipotentes, c'est à dire aux
représentations engendrées par leurs vecteurs fixes sous le sous-groupe d'Iwahori standard $I$. D'après
\cite{vigneras-livre}, le foncteur $V \mapsto V^I$ induit une équivalence de catégories entre la catégorie des
représentations super-unipotentes et la catégorie des modules sur l'algèbre de Hecke-Iwahori. On en déduit alors le
résultat bien connu suivant.

\begin{prop} L'induite normalisée de la représentation triviale du Borel, est semi-simple. Ses facteurs
irréductibles sont en bijection avec les représentations irréductibles du groupe symétrique $\Sigma_d$, la
multiplicité étant alors égale à la dimension de la représentation de $\Sigma_d$ associée.
\end{prop}

\rem la bijection s'explicite comme suit: les représentations irréductibles de $\Sigma_d$ sont en bijection avec les
diagrammes de Young $DY$ associés à une partition $\underline{n}=(d=n_1 \geq n_2 \geq \cdots \geq n_r)$, auxquels on
associe le paramètre de Zelevinski $a(DY)=(1,n_1),\cdots,(1,n_r)$ et donc la représentation irréductible de $GL_d$:
$<a(DY)>$.

\begin{prop} Étant donné un parabolique de Young $\Sigma_{n_1} \times \Sigma_{n_2}$ de $\Sigma_n$ ainsi que
deux diagrammes de Young $DY_1$, $DY_2$
de respectivement $\Sigma_{n_1}$ et $\Sigma_{n_2}$, l'induite parabolique $<a(DY_1)> \times <a(DY_2)>$ est
une somme directe $\bigoplus_{DY} <a(DY)>$ qui porte sur les diagrammes de Young $DY$ qui s'obtiennent à
partir de $DY_1$ en lui adjoignant $\lambda_1$ carrés numérotés $1$ sans en mettre plus d'un dans une même
colonne, puis $\lambda_2$ numérotés $2$ et ainsi de suite où $\lambda_1 \geq \lambda_2 \geq \cdots$ est la
partition de $n_2$ associée à $DY_2$, de sorte qu'en comptant les nouveaux carrés du haut vers le bas et de
droite à gauche, le nombre de $i$ est supérieur ou égal au nombre de $i+1$ pour tout $i$. (cf. \cite{fu-ha}
corollaire 4.39 ainsi que l'appendice).
\end{prop}

\begin{coro} Pour tout $s < l$ et pour tout $0 \leq i \leq s$, la réduction modulo $l$ de
$[\overleftarrow{i},\overrightarrow{s-1-i}]_1$
et $[\overrightarrow{i},\overleftarrow{s-1-i}]_1$ est irréductible.
\end{coro}

\begin{proof} On raisonne par récurrence sur $i$, en partant du fait que pour $i=0$, la réduction modulo $l$
de $[\overleftarrow{s-1}]_1$ et $[\overrightarrow{s-1}]_1$
est irréductible. Supposons donc le résultat acquis jusqu'au rang $i$ et traitons le cas de $i$. L'induite
parabolique $[\overrightarrow{i}]_1 \overrightarrow{\times} [\overleftarrow{s-2-i}]_1$ a pour sous-quotient
$[\overrightarrow{i},\overleftarrow{s-1-i}]_1$ et $[\overrightarrow{i+1},\overleftarrow{s-2-i}]_1$. La
réduction modulo $l$ de cette induite est d'après la proposition précédente la somme directe $<a(DY_1)>
\oplus <a(DY_2)>$ où $DY_1$ (resp. $DY_2$) est le diagramme de Young associé à la partition
$(i+1,1,\cdots,1)$ (resp. $(i+2,1,\cdots,1)$). D'après l'hypothèse de récurrence la réduction modulo $l$ de
$[\overrightarrow{i},\overleftarrow{s-1-i}]_1$ est $<a(DY_1)>$ de sorte que celle de
$[\overrightarrow{i+1},\overleftarrow{s-2-i}]_1$ est $<a(DY_2)>$ et est donc irréductible. Le raisonnement
est identique pour $[\overleftarrow{i},\overrightarrow{s-1-i}]_1$.
\end{proof}

\rem on notera que la réduction modulo $l$ de $[\overleftarrow{1},\overrightarrow{1},\overleftarrow{1}]_1$
n'est pas irréductible mais admet pour sous-quotient irréductible $<a(DY_1)>$ et $<a(DY_2)>$ où $DY_1$
(resp.
$DY_2$) est le diagramme de Young associé à $(2,1,1)$ (resp. $(2,2)$).

\section{Etude de la réduction modulo $l$ de $\st_s(\pi)$}

Dans ce paragraphe $\pi$ désigne une $\bar \Qm_l$-représentation irréductible cuspidale entière de
$Gl_g(K)$ dont on note $\varrho$ la réduction modulo $l$; pour tout $s \geq 1$, on pose $d=sg$.

\subsection{Constituants irréductibles}
\label{para-st-irred}

On rappelle que tous les sous-quotients irréductibles de $\st_s(\pi)$ sont $\varrho$-unipotents;
on se propose alors de les décrire en fonction de leur niveau de cuspidalité au sens de la proposition
\ref{niveau-cusp}.

\begin{prop} \label{prop-super-unique1}
Dans le groupe de Grothendieck des $\bar \Fm_l$-représentations admissibles de $Gl_d(K)$, pour tout $0 \leq t$, la
réduction modulo $l$ de $[\overleftarrow{t}]_{\pi}$ contient une unique représentation $\varrho$-superunipotente que
l'on notera $I_{\underline{0}}([\overleftarrow{t}]_{\pi})$.
\end{prop}

\begin{proof} La preuve procède par récurrence, le cas $t=0$ découlant de l'irréductibilité de $\varrho=r_l(\pi)$ pour
$\pi$ irréductible cuspidale. Pour $t>0$, la réduction modulo $l$ de $J_{P_{tg,g}}([\overleftarrow{t}]_{\pi})=
[\overleftarrow{t-1}]_{\pi\{ 1/2 \} } \otimes [\overleftarrow{0}]_{\pi\{ -t/2 \}}$, d'après l'hypothèse de récurrence,
ne contient qu'une seule représentation $\varrho$-superunipotente, le résultat découle alors  du fait suivant.

Pour une partition $\underline{s}=(s=s_1 \geq s_2 \geq \cdots \geq s_r)$, on note $g.\underline{s}=(d=s_1g \geq s_2 g
\cdots \geq s_r g)$; si $\Pi$ est une représentation $\varrho$-superunipotente, $J_{P_{g.\underline{s}}}(\Pi)$ est non
nulle et contient au moins une représentation $\varrho$-superunipotente.
\end{proof}

\rem d'après le lemme \ref{lem-Zcommute} en utilisant que la réduction modulo $l$ de $[\overrightarrow{t}]_\pi$ est
irréductible et $\varrho$-superunipotente, on en déduit que $I_{\underline{0}}([\overleftarrow{t}]_\pi) \simeq
Z_\varrho([\overrightarrow{t}]_\varrho)$. En outre d'après le corollaire \ref{coro-superunipotent},
$Z_\varrho([\overrightarrow{t}]_\varrho) \simeq Z_1([\overrightarrow{t}]_1) \boxtimes \varrho$ de sorte que
$$I_{\underline{0}}([\overleftarrow{t}]_\pi) \simeq I_{\underline{0}}([\overleftarrow{t}]_1) \boxtimes \varrho \simeq Z_1([\overrightarrow{t}]_1)
\boxtimes \varrho$$ où le calcul de $Z_1([\overrightarrow{t}]_1)$ est donnée dans \cite{ltv}.

Afin d'étudier les autres constituants non $\varrho$-superunipotents, on introduit les ensembles suivants.

\begin{defis} Soient $m_0,\cdots,m_u$ et $s$ des entiers positifs; $\IC(m_0,\cdots,m_u)$ (resp. $\IC_s$) désigne
l'ensemble des $\underline{i}=(i_k)_{k \in \Nm}$ tels que:
\begin{itemize}
\item pour tout $k \geq 0$, $i_k \geq 0$;

\item pour tout $r \in \Nm$, $\sum_{k=r}^\oo (m_k-i_k)l^r \geq 0$ (resp. $s-\sum_{k=0}^\oo i_k m(\varrho)l^k \geq 0$).
\end{itemize}
\end{defis}

\begin{lemm} \label{lem-combi}
Pour $s=m_{-1}+m_0 m(\varrho) + \cdots + m_u m(\varrho)l^u$ avec $0 \leq m_{-1} < m(\varrho)$ et $0 \leq m_0,\cdots,m_u
< l$, on a $\IC(m_0,\cdots,m_u)=\IC_s$ avec pour tout $\underline{i}=(i_k)_{k \in \Nm} \in \IC_s$, $i_k=0$ si $k>u$.
\end{lemm}

\begin{proof}
L'inclusion $\IC_s \subset \IC(m_0,\cdots,m_u)$ est évidente car si $s'=m_{-1}+m(\varrho) \sum_{k=0}^u (m_k-i_k)l^k$
est positif avec $0 \leq m_{-1}< m(\varrho)$ alors $\sum_{k=0}^u (m_k-i_k)l^k \geq 0$. Réciproquement il suffit de
montrer que si $\sum_{k=0}^u (m_k-i_k)l^k \geq 0$ alors $\sum_{k=1}^u (m_k-i_k) l^{k-1} \geq 0$. Or on a
$$\sum_{k=1}^u (m_k-i_k) l^{k-1} \geq \frac{i_0-m_0}{l} > \frac{i_0}{l} -1 >-1$$
et le terme de gauche est un entier qui est donc positif ou nul.
\end{proof}

\begin{defis}
Soit $s=m_{-1}+\sum_{k=0}^u m_k m(\varrho)l^k$ avec $0 \leq m_{-1} < m(\varrho)$ et pour $k=0,\cdots,u$, $0 \leq m_k
<l$:
\begin{itemize}
\item on classe les éléments de $\IC_s=\IC(m_0,\cdots,m_u)$ par ordre lexicographique inverse: $\underline{1}_s=\underline{0}$
est le plus petit élément de $\IC_s$, et le $k$-ème élément de $\IC_s$ par ordre croissant sera noté $\underline{k_s}$.

\item pour $\underline{i} \in \IC_s$, on note $\underline{i}(l)=\sum_{k=0}^\oo i_kl^k$ et
$s_\varrho(\underline{i})=s-m(\varrho)\underline{i}(l)$.
\end{itemize}
\end{defis}

\begin{prop} \label{prop-ss-quotient} La réduction modulo $l$ de
$[\overleftarrow{s-1}]_\pi$ admet $\sharp \IC_s$ constituants irréductibles indexés par les éléments $\underline{i} \in
\IC_s$ de telle sorte que le constituant irréductible $I_{\underline{i}}([\overleftarrow{s-1}]_\pi)$ correspondant est
de niveau de cuspidalité $\underline{i}$ avec les relations suivantes:
\begin{itemize}
\item[(i)] $I_{\underline{0}}([\overleftarrow{s-1}]_\pi)$ est calculé à la proposition \ref{prop-super-unique1};

\item[(ii)] pour tout $\underline{i} \in \IC_s$, on a
$$I_{\underline{i}}([\overleftarrow{s-1}]_\pi) = \Bigl ( I_{\underline{0}}([\overleftarrow{i_0-1}]) \boxtimes \rho_0
\Bigr ) \times \cdots \\ \times \Bigl ( I_{\underline{0}}([\overleftarrow{i_u-1}]) \boxtimes \rho_u \Bigr ) \times
\Bigl ( I_{\underline{0}}([\overleftarrow{s_\varrho(\underline{i})-1}]) \boxtimes \rho_{-1} \Bigr )$$ avec pour
$i=0,\cdots,u$,
$$\rho_i=\st_{m(\varrho)l^i}(\rho\{ \frac{m(\varrho)l^i-s}{2} \} ), \qquad \rho_{-1}=
\varrho \{ \frac{s-s_\varrho(\underline{i})}{2} \};$$

\item[(iii)] pour tout foncteur de Jacquet $J_P$ l'image de
$I_{\underline{i}}([\overleftarrow{s-1}]_\pi)$ est égale à la somme des constituants de niveau de cuspidalité
$\underline{i}$ de $r_l\Bigl ( J_P([\overleftarrow{s-1}]_\pi) \Bigr )$.
\end{itemize}
\end{prop}

\rem pour tout $i=0,\cdots,u$ on a $\rho_i \{ 1 \} \simeq \rho_i$ de sorte que la torsion $\{ n/2 \}$ dans la
définition de $\rho_i$ ne dépend que de la parité de $n$.

\begin{proof}
On raisonne par récurrence sur $s$: les cas $s< m(\varrho)$ découlent directement de \cite{vigneras-livre} car
$r_l([\overleftarrow{s-1}]_\pi)$ est irréductible isomorphe à $\st_s(\varrho)$ qui est $\pi$-superunipotent isomorphe à
$$Z_{\varrho}([\overrightarrow{s-1}]_{\varrho})=Z_1([\overrightarrow{s-1}]_1) \boxtimes \varrho.$$

Supposons donc le résultat acquis pour tout $t<s$ et traitons le cas de $s$. En ce qui concerne les constituants
$\varrho$-superunipotents, le point (i) a été prouvé à la proposition \ref{prop-super-unique1}. Ainsi comme le niveau
de cuspidalité est croissant par foncteur de Jacquet, on en déduit que pour tout parabolique $P$, $J_P \Bigl (
I_{\underline{0}}([\overleftarrow{s-1}]_\pi) \Bigr )$ contient la somme des constituants de niveau $\underline{0}$ de
la réduction modulo $l$ de $J_P \Bigl ( [\overleftarrow{s-1}]_\pi \Bigr )$. En outre
$J_P(I_{\underline{0}}([\overleftarrow{s-1}]_\pi))$ contient au plus un constituant $\varrho$-superunipotent de sorte
que s'il contenait un constituant non $\varrho$-superunipotent, celui-ci serait soit un sous-espace ou un quotient
irréductible de sorte que par réciprocité de Frobenius ou par la deuxième formule d'adjonction,
$I_{\underline{0}}([\overleftarrow{s-1}]_\pi)$ serait un sous-espace ou un quotient irréductible d'une induite d'un
élément non $\varrho$-superunipotent ce qui n'est pas.\footnote{Pour ceux qui voudrait éviter le recours à la deuxième
formule d'adjonction, on peut aussi utiliser un réseau de $[\overleftarrow{t}]_\pi$ et son réseau dual.}

Remarquons ensuite que d'après \cite{vigneras-induced} \S V.3, les représentations
$I_{\underline{i}}([\overleftarrow{s-1}]_\pi)$, définies dans l'énoncé, sont irréductibles. D'après l'hypothèse de
récurrence, les images de $I_{\underline{0}}([\overleftarrow{i_k-1}]) \boxtimes \rho_k$ et $I_{\underline{0}}
([\overleftarrow{s_\varrho(\underline{i})-1}]) \boxtimes \rho_{-1}$ par tout foncteur de Jacquet sont connues, il en
est donc de même pour $I_{\underline{i}}([\overleftarrow{s-1}]_\pi)$. Pour vérifier (iii), il suffit, par transitivité
des foncteurs de Jacquet, de traiter le cas où le parabolique est de la forme $P_{h,d-h}(K)$: ainsi $J_{P_{h,d-h}}
\Bigl ( I_{\underline{i}}([\overleftarrow{s-1}]_\pi) \Bigr )$ est égal à la somme des produits tensoriels:
\begin{multline*}
J_{h_0} \Bigl (I_{\underline{0}}([\overleftarrow{i_0-1}]) \boxtimes \rho_0 \Bigr ) \times \cdots \times J_{h_u} \Bigl
(I_{\underline{0}}([\overleftarrow{i_u-1}]) \boxtimes \rho_u \Bigr ) \times
J_{h_{-1}} \Bigl ( I_{\underline{0}}([\overleftarrow{s_\varrho(\underline{i})-1}]) \boxtimes \rho_{-1} \Bigr ) \\
\otimes J'_{h_0} \Bigl (I_{\underline{0}}([\overleftarrow{i_0-1}]) \boxtimes \rho_0 \Bigr ) \times \cdots \times
J'_{h_u} \Bigl (I_{\underline{0}}([\overleftarrow{i_u-1}]) \boxtimes \rho_u \Bigr ) \times J'_{h_{-1}} \Bigl (
I_{\underline{0}}([\overleftarrow{s_\varrho(\underline{i})-1}] ) \boxtimes \rho_{-1} \Bigr )
\end{multline*}
où la somme porte sur tous les $(h_{-1},h_0,\cdots,h_u)$ tels que:
\begin{itemize}
\item pour tout $k=0,\cdots,u$, on ait $0 \leq h_k \leq gi_km(\varrho)l^k$;

\item $0 \leq h_{-1} \leq s_\varrho(\underline{i}) g$;

\item $h_{-1}+h_0+\cdots+h_u=h$;
\end{itemize}
et où on a posé
$$J_{P_{h_k,gi_km(\rho)l^k-h_k}(K)} \Bigl ( I_{\underline{0}}([\overleftarrow{i_k-1}]) \boxtimes
\rho_k \Bigr )=J_{h_k}  \Bigl (I_{\underline{0}}([\overleftarrow{i_k-1}]) \boxtimes \rho_k \Bigr ) \otimes J'_{h_k}
\Bigl (I_{\underline{0}}([\overleftarrow{i_0-1}]) \boxtimes \rho_0 \Bigr )$$ lesquels sont nuls si pour $k=0,\cdots,u$
(resp. $k=-1$) $h_k$ n'est pas de la forme $gj_k m(\varrho)l^k$ (resp. $h_{-1}=j_{-1}g$) avec $0 \leq j_k \leq i_k$
(resp. $0 \leq j_{-1} \leq s_\varrho(\underline{i})$), et sinon, d'après l'hypothèse de récurrence égaux respectivement
à
$$I_{\underline{0}}([\overleftarrow{j_k-1}]) \boxtimes \rho_k \otimes
I_{\underline{0}}([\overleftarrow{i_k-j_k-1}]) \boxtimes \rho_k$$
$$(\hbox{resp. } I_{\underline{0}}([\overleftarrow{j_{-1}-1}]) \boxtimes \rho_{-1}
\{ \frac{s_\varrho(\underline{i})-j_{-1}}{2} \} \otimes
I_{\underline{0}}[\overleftarrow{s_\varrho(\underline{i})-j_{-1}-1}] \boxtimes \rho_{-1}  \{ -\frac{j_{-1}}{2} \} ).$$
Or d'après l'hypothèse de récurrence, pour $h=tg$,
$$\sum_{t_{-1}+\sum_{k=0}^u t_km(\varrho)l^k=t} I_{\underline{0}}([\overleftarrow{t_0-1}]) \boxtimes
\rho_0 \times \cdots \times I_{\underline{0}}([\overleftarrow{t_{-1}-1}]) \boxtimes \rho_{-1} \{
\frac{s_\varrho(\underline{i})-t_{-1}}{2} \}$$ est égal à la réduction modulo $l$ de $[\overleftarrow{t-1}]_{\pi \{
\frac{s-t}{2} \}}$ tandis que
$$\sum_{t_{-1}+\sum_{k=0}^u t_km(\varrho)l^k=t} I_{\underline{0}}([\overleftarrow{i_0-t_0-1}])
\boxtimes \rho_0 \times \cdots \times I_{\underline{0}}([\overleftarrow{s_\varrho(\underline{i})-t_{-1}-1}]) \boxtimes
\rho_{-1} \{ -\frac{t_{-1}}{2} \}$$ est égal à la réduction modulo $l$ de $[\overleftarrow{s-t-1}]_{\pi \{ -\frac{t}{2}
\} }$. En résumé tous les constituants de $J_{P_{tg,(s-t)g}} \Bigl ( I_{\underline{i}}([\overleftarrow{s-1}]_\pi) \Bigr
)$ sont de niveau de cuspidalité $\underline{i}$ et correspondent aux constituants de niveau de cuspidalité
$\underline{i}$ de la réduction modulo $l$ de $[\overleftarrow{t-1}]_{\pi \{ (s-t)/2 \} } \otimes
[\overleftarrow{s-t-1}]_{\pi \{ -t/2 \} }$ ce qui prouve le point (iii).

Supposons désormais avoir montré par récurrence que pour tout $\underline{i}<\underline{k_s}$,
$I_{\underline{i}}([\overleftarrow{s-1}]_\pi)$ tel qu'il est défini ci-dessus, est un constituant de la réduction
modulo $l$ de $[\overleftarrow{s-1}]_\pi$. Les calculs ci-dessus montrent que dans le groupe de Grothendieck,
$\Psi_k:=r_l([\overleftarrow{s-1}]_\pi) - \sum_{\underline{i} < \underline{k_s}}
I_{\underline{i}}([\overleftarrow{s-1}]_\pi)$ est un élément effectif du groupe de Grothendieck tel que son image par
tout foncteur de Jacquet propre, est de niveau de cuspidalité supérieur ou égal à $\underline{k_s}$. Soit alors $-1
\leq j \leq u$ minimal tel que $\rho_j$, à torsion près, appartienne au support cuspidal de $I_{\underline{k_s}}
([\overleftarrow{s-1}]_\pi)$; pour $j \geq 0$ (resp. $j=-1$) on note $P_j$ le parabolique standard de groupe de Levi
$GL_{(s-m(\varrho)l^j)g} \times GL_{m(\varrho)l^j g}$ (resp. $GL_{(s-1)g} \times GL_g$) de sorte que $J_{P_j}(\Psi_k)$
est égal à
$$I_{\underline{k_s}-(0,\cdots,0,1,0,\cdots)}([\overleftarrow{s-1-m(\varrho)l^j}]_\pi) \otimes \rho_j \quad
(\hbox{resp. } I_{\underline{k_s}}([\overleftarrow{s-2}]_{\pi\{ \frac{1}{2} \} }) \otimes \rho_{-1} \{ \frac{1-s}{2}
\})$$ qui est son unique constituant de niveau de cuspidalité $\underline{k_s}$. On en déduit alors que $\Psi_k$
contient un unique constituant $\pi_k$ de niveau de cuspidalité $\underline{k_s}$ qui, par réciprocité de Frobenius,
est un sous-espace de l'induite parabolique
$$I_{\underline{k_s}-(0,\cdots,0,1,0,\cdots)}([\overleftarrow{s-1-m(\varrho)l^j}]_\pi) \times \rho_j \quad
(\hbox{resp. } I_{\underline{k_s}}([\overleftarrow{s-2}]_{\pi \{ 1/2 \} }) \times \rho_{-1} \{ \frac{1-s}{2} \}).$$ En
utilisant que l'image de $\Psi_k$ par le foncteur de Jacquet associé au parabolique de Levi $GL_{(s-k_j m(\varrho)
l^j)g} \times GL_{k_jm(\varrho)l^j}$ (resp. $GL_{(s-s_\varrho(\underline{i}))g} \times GL_{s_\varrho(\underline{i})g}$)
admet un unique constituant tel que son deuxième facteur est de niveau de cuspidalité $(0,\cdots,0,k_j,0,\cdots)$, on
obtient que $\pi_k$ est nécessairement égal à $I_{\underline{k_s}}([\overleftarrow{s-1}]_\pi)$.

Au final $r_l([\overleftarrow{s-1}]_\pi) - \sum_{\underline{i} \in \IC_s} I_{\underline{i}}([\overleftarrow{s-1}]_\pi)$
est un élément effectif du groupe de Grothendieck dont toutes les images par les foncteurs de Jacquet propres sont
nulles. Le résultat découle alors du fait que $r_l([\overleftarrow{s-1}]_\pi)$ ne contient aucun constituant cuspidal
pour $s$ qui n'est pas de la forme $m(\varrho)l^i$ et que sinon $\st_{m(\varrho)l^i}(\varrho)$ en est l'unique
constituant cuspidal.
\end{proof}

\rem contrairement au cas général, tous les constituants irréductibles $\psi$ de $r_l([\overleftarrow{s-1}]_\pi)$ sont
tels que leur niveau de cuspidalité reste constant par foncteur de Jacquet, i.e. pour tout parabolique $P$, le niveau
de cuspidalité d'un constituant irréductible de $J_P(\psi)$ est égal à celui de $\psi$.

\begin{defi} Pour $k \geq 0$ (resp. $k=-1$), on dit que qu'une représentation de $GL_{sg}(K)$ $\varrho$-unipotente est de $\varrho_k$-niveau de
cuspidalité nul si son niveau de cuspidalité $\underline{i}=(i_k)_{k \in \Nm}$ est tel que $i_k=0$ (resp.
$s_\varrho(\underline{i})=0$).
\end{defi}

\subsection{Constructions des réseaux d'induction}
\label{para-st-reseau}

Pour $\pi$ une représentation irréductible entière cuspidale de $GL_g(K)$, on a vu que sa réduction modulo $l$,
$\varrho$ était irréductible de sorte qu'à isomorphismes près, $[\overleftarrow{0}]_\pi$ possède un unique réseau
stable. Ainsi étant donné un réseau de $[\overleftarrow{t-1}]_\pi$, la surjection $[\overleftarrow{t-1}]_\pi
\overrightarrow{\times} [\overleftarrow{0}]_\pi \twoheadrightarrow [\overleftarrow{t}]_\pi$ induit un réseau de
$[\overleftarrow{t}]_\pi$ de sorte que par récurrence on dispose d'un réseau particulier $RI_{\bar
\Zm_l}([\overleftarrow{t-1}]_\pi)$ que l'on qualifie \textit{réseau d'induction}, des $[\overleftarrow{t-1}]_\pi$. La
motivation principale pour l'étude de ces réseaux d'induction est qu'ils apparaissent dans la cohomologie des espaces
de Deligne-Carayol.

\begin{prop} \label{prop-reseau}
Il existe un réseau de $[\overleftarrow{s-1}]_\pi$, unique à isomorphisme près, que l'on notera $RI_{\bar
\Zm_l}([\overleftarrow{s-1}]_\pi)$, tel que dans son graphe $\Gamma_{RI_{\bar \Zm_l}([\overleftarrow{s-1}]_\pi)}$, les
seules flèches sont celles qui relient $I_{\underline{i}}([\overleftarrow{s-1}]_\pi)$ à
$I_{\underline{i'}}([\overleftarrow{s-1}]_\pi)$, où $\underline{k_s}=\underline{i} <
\underline{(k-1)_s}=\underline{i'}$ sont des éléments consécutifs de $\IC_s$. En outre $RI_{\bar
\Zm_l}([\overleftarrow{s-1}]_\pi)$ vérifie les propriétés suivantes:
\begin{itemize}
\item[(i)] si $\Lambda$ est un réseau stable de $[\overleftarrow{s}]_\pi$ tel que la surjection
$$[\overleftarrow{s-1}]_\pi \overrightarrow{\times} [\overleftarrow{0}]_\pi \twoheadrightarrow
[\overleftarrow{s}]_\pi$$ induise un morphisme surjectif de $RI_{\bar \Zm_l}([\overleftarrow{s-1}]_\pi)
\overrightarrow{\times} RI_{\bar \Zm_l}([\overleftarrow{0}]_\pi)$ sur $\Lambda$ alors $\Lambda \simeq RI_{\bar
\Zm_l}([\overleftarrow{s}]_\pi)$;

\item[(ii)] pour tout $1 \leq t \leq s$, l'image de $RI_{\bar \Zm_l}([\overleftarrow{s-1}]_\pi)$
par le foncteur de Jacquet $J_{P_{tg,(s-t)g}}$ est isomorphe à $RI_{\bar \Zm_l}([\overleftarrow{t-1}]_{\pi \{
\frac{t-s}{2} \} }) \otimes RI_{\bar \Zm_l}([\overleftarrow{s-t-1}]_{\pi \{ \frac{t}{2} \} })$.
\end{itemize}
\end{prop}

\begin{proof}
On raisonne par récurrence sur $s$. Pour $1 \leq s < m(\varrho)$ on rappelle que $r_l([\overleftarrow{s-1}]_\pi)$ est
irréductible, il n'y a donc qu'un seul réseau; il s'agit alors simplement de vérifier (i) dans le cas $s=m(\varrho)-1$.
Dans le groupe de Grothendieck, $r_l([\overleftarrow{s}]_\pi)$ est de longueur $2$ avec un constituant
$\varrho$-superunipotent et l'autre cuspidal qui ne peut donc pas être un quotient de $[\overleftarrow{s-1}]_{\varrho}
\overrightarrow{\times} [\overleftarrow{0}]_{\varrho}$: ainsi le réseau image de $RI_{\bar
\Zm_l}([\overleftarrow{s-1}]_\pi) \overrightarrow{\times} RI_{\bar \Zm_l}([\overleftarrow{0}]_\pi)$ dans
$[\overleftarrow{s}]_\pi$ définit une extension non triviale de $[\overleftarrow{s}]_{\varrho}$ par
$I_{\underline{0}}([\overleftarrow{s}]_\pi)$.

Supposons avoir construit pour tout $1 \leq t<s$, un réseau $RI_{\bar \Zm_l}([\overleftarrow{t-1}]_\pi)$ de
$[\overleftarrow{t-1}]_\pi$ dont le graphe d'extension est comme dans l'énoncé et vérifiant (i) (resp. (ii)) pour tout
$t<s-1$ (resp. $t<s$). On définit alors $RI_{\bar \Zm_l}([\overleftarrow{s-1}]_\pi)$ comme le réseau image de $RI_{\bar
\Zm_l}([\overleftarrow{s-2}]_\pi) \overrightarrow{\times} RI_{\bar \Zm_l}([\overleftarrow{0}]_\pi)$ via l'application
surjective $[\overleftarrow{s-2}]_\pi \overrightarrow{\times} [\overleftarrow{0}]_ \pi \twoheadrightarrow
[\overleftarrow{s-1}]_\pi$ et on note
$$RI_{\bar \Fm_l}([\overleftarrow{s-1}]_\pi):= RI_{\bar \Zm_l}([\overleftarrow{s-1}]_\pi) \otimes_{\bar
\Zm_l}
\bar \Fm_l.$$

Pour tout $1 \leq t < s$, le foncteur de Jacquet $J_{P_{tg,(s-t)g}}$ étant exact à droite, on en déduit une application
surjective
$$J_{P_{tg,(s-t)g}}\Bigl ( RI_{\bar \Zm_l}([\overleftarrow{s-2}]_\pi) \overrightarrow{\times}
RI_{\bar \Zm_l}([\overleftarrow{0}]_\pi) \Bigr ) \twoheadrightarrow J_{P_{tg,(s-t)g}}\Bigl ( RI_{\bar
\Zm_l}([\overleftarrow{s-1}]_\pi) \Bigr ).$$ D'après \cite{ze} et (ii), le membre de gauche possède une filtration de
longueur $2$ dont les gradués sont
\begin{itemize}
\item $RI_{\bar \Zm_l}([\overleftarrow{t-1}]_{\pi \{ \frac{s-t-2}{2} \} }) \otimes
RI_{\bar \Zm_l}([\overleftarrow{s-t-2}]_{\pi \{ -\frac{t+1}{2} \} }) \times RI_{\bar \Zm_l}([\overleftarrow{0}]_{\pi\{
\frac{s-1}{2} \} })$,
\item $RI_{\bar \Zm_l}([\overleftarrow{t-2}]_{\pi \{ \frac{s-t}{2} \} }) \overrightarrow{\times}
RI_{\bar \Zm_l}([\overleftarrow{0}]_{\pi \{ \frac{s-t}{2} \} }) \otimes RI_{\bar \Zm_l}([\overleftarrow{s-t-1}]_{\pi\{
-t/2 \} })$.
\end{itemize}
D'après l'hypothèse de récurrence et la propriété (i), la réduction modulo $l$ de chacun de ces gradués possède un
unique quotient irréductible à savoir $I_{\underline{0}} ([\overleftarrow{t-1}]_\pi) \otimes I_{\underline{0}}
([\overleftarrow{s-t-1}]_\pi)$ de sorte que $J_{P_{tg,(s-t)g}}\Bigl ( RI_{\bar \Zm_l}([\overleftarrow{s-1}]_\pi) \Bigr
)$ est un quotient de l'une des deux représentations induites suivantes, le choix dépendant à priori de $t$:
\begin{equation} \label{eq-choix}
J_{P_{tg,(s-t)g}}\Bigl ( RI_{\bar \Zm_l}([\overleftarrow{s-1}]_\pi) \Bigr ) \twoheadleftarrow \left \{ \begin{array}{l}
RI_{\bar \Zm_l}([\overleftarrow{t-1}]_{\pi \{ \frac{s-t-2}{2} \} }) \otimes RI_{\bar
\Zm_l}([\overleftarrow{s-t-2}]_{\pi \{ -\frac{t+1}{2} \} }) \times
RI_{\bar \Zm_l}([\overleftarrow{0}]_{\pi\{ \frac{s-1}{2} \} }) \\
RI_{\bar \Zm_l}([\overleftarrow{t-2}]_{\pi \{ \frac{s-t}{2} \} }) \overrightarrow{\times} RI_{\bar
\Zm_l}([\overleftarrow{0}]_{\pi \{ \frac{s-t}{2} \} }) \otimes RI_{\bar \Zm_l}([\overleftarrow{s-t-1}]_{\pi\{ -t/2 \}
})
\end{array} \right .
\end{equation}
On en déduit alors que $I_{\underline{0}}([\overleftarrow{s-1}]_\pi)$ est l'unique quotient irréductible de $RI_{\bar
\Fm_l}([\overleftarrow{s-1}]_\pi)$. En effet $[\overleftarrow{s-1}]_{\varrho}$ est l'unique constituant cuspidal de
$r_l([\overleftarrow{s-1}]_\pi)$ pour $s$ de la forme $m(\varrho)l^u$ et il n'est alors pas un quotient de $RI_{\bar
\Fm_l}([\overleftarrow{s-1}]_\pi) \overrightarrow{\times} RI_{\bar \Fm_l}([\overleftarrow{0}]_\pi)$. Soit alors
$I_{\underline{i}}([\overleftarrow{s-1}]_\pi)$ un quotient de $RI_{\bar \Fm_l}([\overleftarrow{s-1}]_\pi)$ de sorte
qu'il existe $1 \leq t <s$ tel que $J_{P_{tg,(s-t)g}} \Bigl ( I_{\underline{i}}([\overleftarrow{s-1}]_\pi) \Bigr )$
soit un quotient de $J_{P_{tg,(s-t)g}}\Bigl ( RI_{\bar \Fm_l}([\overleftarrow{s-1}]_\pi) \Bigr )$; d'après ce qui
précède $J_{P_{tg,(s-t)g}} \Bigl ( I_{\underline{i}}([\overleftarrow{s-1}]_\pi) \Bigr )$ est alors
$\varrho$-superunipotent ce qui impose $\underline{i}=\underline{0}$.

Notons alors
$$V_1(s)= \ker \Bigl ( RI_{\bar \Fm_l}([\overleftarrow{s-1}]_\pi) \twoheadrightarrow
I_{\underline{0}}([\overleftarrow{s-1}]_\pi) \Bigr ),$$ et supposons avoir montré par récurrence l'unicité d'une suite
$$V_{k-1}(s) \varsubsetneq V_{k-2}(s) \varsubsetneq \cdots \varsubsetneq V_1(s) \varsubsetneq
V_{0}(s)=RI_{\bar \Fm_l}([\overleftarrow{s-1}]_\pi)$$ telle que pour tout $0 \leq i \leq k-1$, $V_i(s)$ est l'unique
sous-espace stable strict de $V_{i-1}(s)$ tel que $V_{i-1}(s)/V_{i}(s)$ soit irréductible et qu'alors ce quotient est
isomorphe à $I_{\underline{i_s}}([\overleftarrow{s-1}]_\pi)$. Montrons alors qu'il existe un unique sous-espace stable
strict $V_k(s)$ de $V_{k-1}(s)$ telle que $V_{k-1}(s)/V_k(s)$ soit irréductible et qu'alors ce quotient est isomorphe à
$I_{\underline{k_s}}([\overleftarrow{s-1}]_\pi)$. Soit $k' \geq 1$ tel que $\underline{(k-1)_s} < \underline{k'_{s-1}}
\leq \underline{k_s}$ de sorte que la composée
\begin{multline*}
V_{k'}(s-1) \overrightarrow{\times} RI_{\bar \Fm_l}([\overleftarrow{0}]_\pi) \hookrightarrow RI_{\bar
\Fm_l}([\overleftarrow{s-2}]_\pi) \overrightarrow{\times} RI_{\bar \Fm_l}([\overleftarrow{0}]_\pi) \\
\twoheadrightarrow RI_{\bar \Fm_l}([\overleftarrow{s-1}]_\pi) \twoheadrightarrow RI_{\bar
\Fm_l}([\overleftarrow{s-1}]_\pi)/ V_{k-1}(s)
\end{multline*}
est nulle car par construction tous les constituants de $V_{k'}(s-1)$ sont de niveau de cuspidalité supérieur ou égal à
$\underline{k'_{s-1}}$ et donc aussi ceux de $V_{k'}(s-1) \overrightarrow{\times} RI_{\bar
\Fm_l}([\overleftarrow{0}]_\pi)$. On en déduit donc une application
$$V_{k'}(s-1) \overrightarrow{\times} RI_{\bar \Fm_l}([\overleftarrow{0}]_\pi) \rightarrow
V_{k-1}(s).$$ Si $r_l([\overleftarrow{s-1}]_\pi)$ contient une cuspidale alors $s$ est de la forme $m(\varrho)l^u$ et
il s'agit de $[\overleftarrow{s-1}]_{\varrho}$. On note aussi que $[\overleftarrow{s-1}]_{\varrho}$ est un constituant
de multiplicité $1$ de la réduction modulo $l$ de $[\overleftarrow{s-2}]_\pi \overrightarrow{\times}
[\overleftarrow{0}]_\pi$; c'est aussi un constituant de $V_{k'}(s-1) \overrightarrow{\times} RI_{\bar
\Fm_l}([\overleftarrow{0}]_\pi)$. Ainsi si la représentation cuspidale $[\overleftarrow{s-1}]_{\varrho}$ était un
quotient de $V_{k-1}(s)$, du fait de la surjectivité de $RI_{\bar \Fm_l}([\overleftarrow{s-2}]_\pi)
\overrightarrow{\times} RI_{\bar \Fm_l}([\overleftarrow{0}]_\pi) \twoheadrightarrow RI_{\bar
\Fm_l}([\overleftarrow{s-1}]_\pi)$, il serait aussi un quotient de $V_{k'}(s-1) \overrightarrow{\times} RI_{\bar
\Fm_l}([\overleftarrow{0}]_\pi)$ ce qui ne se peut pas.

Soit alors $I_{\underline{i}}([\overleftarrow{s-1}]_\pi)$ un quotient de $V_{k-1}(s)$ et $P_{tg,(s-t)g}$ un parabolique
tel que
$$J_{P_{tg,(s-t)g}} \Bigl ( I_{\underline{i}}([\overleftarrow{s-1}]_\pi) \Bigr ) \neq 0.$$
Par construction on a
$\underline{i} \geq \underline{k_s}$ et dans le groupe de Grothendieck $J_{P_{tg,(s-t)g}} \Bigl ( V_{k-1}(s) \Bigr )$
est la somme des constituants de niveau de cuspidalité supérieurs ou égaux à $\underline{k_s}$ de
$r_l([\overleftarrow{t-1}]_\pi) \otimes r_l([\overleftarrow{s-t-1}]_\pi)$. En outre d'après l'hypothèse de récurrence
appliquée à (\ref{eq-choix}) montre que $J_{P_{tg,(s-t)g}} \Bigl ( V_{k-1}(s) \Bigr )$ a un unique quotient
irréductible, celui de plus bas niveau de cuspidalité. Il suffit alors de montrer que si
$I_{\underline{i'}}([\overleftarrow{s-1}]_\pi)$ est un quotient de $V_{k-1}(s)$ alors il existe un parabolique
$P_{tg,(s-t)g}(K)$ tel que $J_{N_{tg,(s-t)g}} \Bigl ( I_{\underline{i}}([\overleftarrow{s-1}]_\pi) \Bigr ) \neq 0$ et
$J_{N_{tg,(s-t)g}} \Bigl ( I_{\underline{i'}}([\overleftarrow{s-1}]_\pi) \Bigr ) \neq 0$.

Posons $i_{-1}=s-m(\varrho) \sum_{k=0}^u i_kl^k$ et $i'_{-1}=s-m(\varrho) \sum_{k=0}^u i'_kl^k$; s'il existe $0 \leq k$
(resp. $k=-1$) tel que $i_ki'_k \neq 0$ alors $t=m(\varrho)l^k$ (resp. $t=1$) convient: en effet on a alors
\begin{multline*}
J_{P_{tg,(s-t)g}}(I_{\underline{i}}([\overleftarrow{s-1}]_\pi))= \rho_k \otimes
I_{\underline{i}-(0,\cdots,1,0,\cdots)}([\overleftarrow{s-t-1}]_\pi) \\
(\hbox{resp. } \rho_{-1} \{ \frac{s-1}{2} \} \otimes I_{\underline{i}}([\overleftarrow{s-2}]_{\pi \{ -\frac{1}{2} \}
}))
\end{multline*}
et de même pour $I_{\underline{i'}}([\overleftarrow{s-1}]_\pi)$. Supposons donc que pour tout $k=-1,\cdots,u$, on ait
$i_ki'_k=0$ et notons $-1 \leq k_0 \leq u$ le plus petit indice $k$ tel que $i_k \neq 0$: soit alors $r$ tel que
$i'_{k_0+r}$ soit le premier des $i'_k$ non nul. Quitte à échanger le rôle de $\underline{i}$ et $\underline{i'}$, on
peut supposer $r >0$ de sorte que $t=m(\varrho)l^{k_0+r}$ convient. En effet dans le cas où $k_0 \geq 0$ alors $l^r |
\sum_{w=0}^{r-1} i_{k_0+w}l^w$ et donc $\sum_{w=0}^{r-1} i_{k_0+w}l^w \geq l^r$: soit d'après le lemme \ref{lem-combi},
$(j_{k_0},\cdots,j_{k_0+r-1})$ tels que:
\begin{itemize}
\item pour tout $w=0,\cdots,r-1$, $0 \leq j_{k_0+w} \leq i_{k_0+w}$ et

\item $\sum_{w=0}^{r-1} j_{k_0+w}l^w=l^r$,
\end{itemize}
de sorte que $J_{P_{tg,(s-t)g}}(I_{\underline{i}}([\overleftarrow{s-1}]_\pi))$ admet comme constituant
\begin{multline*}
I_{\underline{0}}([\overleftarrow{j_{k_0}-1}]) \boxtimes \rho_{k_0} \times \cdots \times
I_{\underline{0}}([\overleftarrow{j_{k_0+r-1}-1}]) \boxtimes \rho_{k_0+r-1} \otimes \\
I_{\underline{0}}([\overleftarrow{i_{k_0}-j_{k_0}-1}]) \boxtimes \rho_{k_0} \times \cdots \times
I_{\underline{0}}([\overleftarrow{i_{k_0+r-1}-j_{k_0+r-1}-1}]) \boxtimes \rho_{k_0+r-1} \\ \times
I_{\underline{0}}([\overleftarrow{i_{k_0+r}-1}]) \boxtimes \rho_{k_0+r} \times \cdots \times
I_{\underline{0}}([\overleftarrow{i_u-1}]) \boxtimes \rho_{u}
\end{multline*}
 alors que $J_{P_{tg,(s-t)g}}(I_{\underline{i}}([\overleftarrow{s-1}]_\pi))$ admet $\rho_{k_0+r}
\otimes I_{(i'_0,\cdots,i'_{k_0+r}-1,i'_{k_0+r+1},\cdots)}([\overleftarrow{s-t-1}]_\pi)$.

Si $k_0=-1$ alors $i_{-1}+m(\varrho)\sum_{k=0}^{r-1} i_k l^k$ est divisible par $m(\varrho)l^r$; le reste du
raisonnement est ensuite strictement identique au cas $k_0 \geq 0$ en prenant soin de préciser les torsions sur
$\rho_{-1}$.
\end{proof}

\begin{coro} \label{coro-st-induite}
Soit $\pi$ une représentation irréductible cuspidale entière de $GL_g(K)$ et soient deux réseaux
stable $R_{\bar \Zm_l}([\overleftarrow{t}]_\pi)$ et $R_{\bar \Zm_l}([\overleftarrow{a}]_\pi)$ telle que l'on
ait une surjection
$$R_{\bar \Zm_l}([\overleftarrow{t}]_\pi) \overrightarrow{\times }R_{\bar \Zm_l}([\overleftarrow{a}]_\pi)
\twoheadrightarrow RI_{\bar \Zm_l}([\overleftarrow{t+a+1}]_\pi)$$ Alors les deux réseaux en question sont isomorphes
respectivement à $RI_{\bar \Zm_l}([\overleftarrow{t}]_\pi)$ et $RI_{\bar \Zm_l}([\overleftarrow{a}]_\pi)$.
\end{coro}

\begin{proof}
Par application du foncteur de Jacquet $J_{P_{(t+1)g,(l+a+2)g}(K)}$, on en déduit une application surjective
$$R_{\bar \Zm_l}([\overleftarrow{t}]_{\pi\{ -a/2 \}}) \otimes R_{\bar \Zm_l}([\overleftarrow{a}]_{\pi\{ t/2 \} })
\twoheadrightarrow RI_{\bar \Zm_l}([\overleftarrow{t}]_{\pi\{ -a/2 \} }) \otimes RI_{\bar \Zm_l}
([\overleftarrow{a}]_{\pi\{t/2\} })$$ d'où le résultat.
\end{proof}

\section{Etude de la réduction modulo $l$ de $LT_t(\pi,s)$}
\label{para-red2}

Pour $\pi$ une $\bar \Qm_l$-représentation irréductible entière cuspidale de $GL_g(K)$, on se propose dans ce
paragraphe de reprendre les résultats précédents pour $LT_t(\pi,s):=[\overleftarrow{t-1},\overrightarrow{s-t}]_\pi$. On
note comme précédemment $\varrho=r_l(\pi)$ et $d=sg$.

\subsection{Constituants irréductibles}

Commençons comme précédemment par étudier les constituants $\varrho$-superunipotents.

\begin{prop} \label{prop-super-unique2}
Dans le groupe de Grothendieck des $\bar \Fm_l$-représentations admissibles de $Gl_d(K)$, pour tout $0 \leq t \leq a$,
la réduction modulo $l$ de $[\overleftarrow{t},\overrightarrow{a-t}]_{\pi}$ contient une unique représentation
$\varrho$-superunipotente que l'on notera $I_{\underline{0}}([\overleftarrow{t},\overrightarrow{a-t}]_{\pi})$. Par
ailleurs dans le groupe de Grothendieck des $\bar \Fm_l$-représentations admissibles de $Gl_d(K)$, l'image de
$I_{\underline{0}}([\overleftarrow{t},\overrightarrow{a-t}]_{\pi})$ par le foncteur de Jacquet $J_{P_{(t+1)g,(a-t)g}}$
(resp. $J_{P_{(a-t+1)g,tg}}$, resp. $J_{P_{at,g}}$) contient
\begin{multline*}
I_{\underline{0}}([\overleftarrow{t}]_{\pi\{ \frac{t-a}{2}\} }) \otimes [\overrightarrow{a-t-1}]_{\varrho\{ \frac{t+1}{2}\} }), \\
(\hbox{resp. }[\overrightarrow{a-t}]_{\varrho \{ -t/2 \} } \otimes
I_{\underline{0}}([\overleftarrow{t-1}]_{\pi\{ \frac{a-t+1}{2} \} }) \\
\hbox{resp. } I_{\underline{0}}([\overleftarrow{t},\overrightarrow{a-t-1}]_{\pi\{-1/2 \} }) \otimes
[\overrightarrow{0}]_{r_l(\pi\{ a/2 \} )} + \\
I_{\underline{0}}([\overleftarrow{t-1},\overrightarrow{a-t}]_{\pi\{ 1/2 \} }) \otimes [\overleftarrow{0}]_{\varrho \{
-a/2 \} })
\end{multline*}
\end{prop}

\begin{proof} Le cas $a=t$ a été traité à la proposition \ref{prop-super-unique1}. Par ailleurs rappelons que le
cas $t=0$ découle du fait que d'après \cite{vigneras-livre} V.9.1 (c), $r_l([\overrightarrow{a}]_\pi)$ est irréductible
et $\varrho$-superunipotent.\footnote{Il est aussi possible de raisonner par récurrence: pour $a=0$, le résultat se
déduit du fait que $\varrho$ est irréductible et pour $a>0$, la réduction modulo $l$ de
$J_{P_{ag,g}}([\overrightarrow{a}]_{\pi})= [\overrightarrow{a-1}]_{\pi\{-1/2\} } \otimes [\overrightarrow{0}]_{\pi\{
a/2 \} }$ contient, d'après l'hypothèse de récurrence un unique constituant $\varrho$-superunipotent, ce qui implique
le résultat.}

On raisonne par récurrence en utilisant la remarque suivante de la proposition \ref{prop-super-unique1}. Soit
$\underline{s}=(s=s_1 \geq s_2 \geq \cdots \geq s_r)$ une partition de $s$ et  $g.\underline{s}=(d=s_1g \geq s_2 g
\cdots \geq s_r g)$: si $\pi$ est une représentation $\varrho$-superunipotente alors $J_{P_{g.\underline{s}}}(\pi)$ est
non nulle et contient au moins une représentation $\varrho$-superunipotente. Par application du foncteur de Jacquet
$J_{P_{ag,g}}$, on obtient dans le groupe de Grothendieck
$$J_{P_{ag,g}}([\overleftarrow{t},\overrightarrow{a-t}]_\pi)=[\overleftarrow{t},\overrightarrow{a-t-1}]_{\pi
\{ -1/2 \} } \otimes \pi \{ a/2 \} + [\overleftarrow{t-1},\overrightarrow{a-t}]_{\pi \{ 1/2 \} } \otimes \pi \{ -a/2
\}$$ de sorte que, d'après l'hypothèse de récurrence et la remarque ci-dessus, la réduction modulo $l$ de
$[\overleftarrow{t},\overrightarrow{a-t}]_\pi$ contient un ou deux constituants $\varrho$-superunipotent.

\begin{lemm} \label{lem-reseau0}
Les réseaux de $[\overleftarrow{t},\overrightarrow{a-t}]_\pi$ induits par
$$RI_{\bar \Zm_l}([\overleftarrow{t}]_\pi) \overrightarrow{\times} [\overrightarrow{a-t-1}]_\pi \hbox{ et }[\overrightarrow{a-t}]_\pi \overleftarrow{\times} RI_{\bar
\Zm_l}([\overleftarrow{t-1}]_\pi)$$ sont isomorphes.
\end{lemm}

\begin{proof}
On raisonne par récurrence sur $a-t$; le cas $a-t=0$ découle de la caractérisation des réseaux $RI_{\bar
\Zm_l}([\overleftarrow{t}]_\pi)$ donnée à la proposition \ref{prop-reseau}. Supposons donc le résultat acquis jusqu'au
rang $a-t-1$ et traitons le cas de $a-t$. On rappelle que sur $\bar \Qm_l$, on a
$$\hom_{\bar \Qm_l[GL_{(a+1)g}(K)]} \Bigl ( [\overleftarrow{t}]_\pi \overrightarrow{\times} [\overrightarrow{a-t-1}]_\pi,
[\overrightarrow{a-t}]_\pi \overleftarrow{\times} [\overleftarrow{t-1}]_\pi \Bigr ) \simeq \bar \Qm_l$$ toute flèche
non nulle induisant un isomorphisme de $[\overleftarrow{t},\overrightarrow{a-t}]_\pi$ sur lui-même. Au niveau des
réseaux, d'après la réciprocité de Frobenius, on a
\begin{multline*}
\hom_{\bar \Zm_l[GL_{(a+1)g}(K)]} \Bigl ( RI_{\bar \Zm_l}([\overleftarrow{t}]_\pi) \overrightarrow{\times}
[\overrightarrow{a-t-1}]_\pi, [\overrightarrow{a-t}]_\pi \overleftarrow{\times} RI_{\bar
\Zm_l}([\overleftarrow{t-1}]_\pi) \Bigr ) \simeq \\
\hom_{\bar \Zm_l[GL_{(a+1)g}(K)]} \Bigl ( J_{P_{(a-t+1)g,tg}} (
RI_{\bar \Zm_l} ([\overleftarrow{t}]_\pi) \overrightarrow{\times} [\overrightarrow{a-t-1}]_\pi) ,
[\overrightarrow{a-t}]_{\pi \{ t/2 \} } \otimes RI_{\bar \Zm_l} ([\overleftarrow{t-1}]_{\pi \{ (t-a-1)/2 \} } \Bigr )
\end{multline*}
où par égalité des supports cuspidaux, ce dernier espace est
$$\hom_{\bar \Zm_l[GL_{(a+1)g}(K)]} \Bigl ( RI_{\bar \Zm_l}([\overleftarrow{0}]_{\pi \{ t/2 \} }) \overrightarrow{\times}
[\overrightarrow{a-t-1}]_{\pi \{ t/2 \} } \otimes RI_{\bar \Zm_l}([\overleftarrow{t-1}]_{\pi \{ (t-a-1)/2 \} }) \Bigr )
\simeq \bar \Zm_l$$ car d'après \ref{prop-reseau}, $J_{P_{g,tg}} ( RI_{\bar \Zm_l} ([\overleftarrow{t}]_\pi)) \simeq
RI_{\bar \Zm_l}([\overleftarrow{0}]_{\pi \{ t/2 \} } \otimes RI_{\bar \Zm_l} ([\overleftarrow{t-1}]_{\pi \{ -1/2 \} }
)$, d'où le résultat.
\end{proof}

Pour le réseau du lemme précédent, notons $\rho$ une sous-représentation irréductible de
$r_l([\overleftarrow{t},\overrightarrow{a-t}]_\pi)$ de sorte que $\rho$ est aussi une sous-représentation irréductible
des réductions modulo $l$ de $RI_{\bar \Zm_l}([\overleftarrow{t}]_\pi) \overrightarrow{\times}
[\overrightarrow{a-t-1}]_\pi$ et $[\overrightarrow{a-t}]_\pi \overleftarrow{\times} RI_{\bar
\Zm_l}([\overleftarrow{t-1}]_\pi)$. Par réciprocité de Frobenius, on en déduit alors que l'image de $\rho$ par
$J_{P_{(t+1)g,(a-t)g}}$ (resp. $J_{P_{(a-t+1)g,tg}}$) contient $I_{\underline{0}}([\overleftarrow{t}]_{\pi \{ (t-a)/2
\} }) \otimes [\overrightarrow{a-t-1}]_{\varrho \{ (t+1)/2 \} }$ (resp. $[\overrightarrow{a-t}]_{\varrho \{ t/2 \} }
\otimes I_{\underline{0}}([\overleftarrow{t-1}]_{\pi \{ (t-1-a)/2 \} })$) de sorte que $J_{P_{g,\cdots,g}}(\rho)$
contient $? \otimes \cdots \otimes ? \otimes \varrho \{ a/2 \}$ (resp. $? \otimes \cdots \otimes ? \otimes \varrho \{
-a/2 \}$). Ainsi $J_{P_{ag,g}}(\rho)$ contient forcément $I_{\underline{0}} \Bigl (
[\overleftarrow{t},\overrightarrow{a-t-1}]_{\pi \{ -1/2 \} } \Bigr ) \otimes \pi \{ a/2 \}$ (resp. $I_{\underline{0}}
\Bigl ( [\overleftarrow{t-1},\overrightarrow{a-t}]_{\pi \{ 1/2 \} } \Bigr ) \otimes \pi \{ -a/2 \}$) ce qui prouve
l'unicité du constituant $\varrho$-suerunipotent de la réduction modulo $l$ de
$[\overleftarrow{t},\overrightarrow{a-t}]_\pi$.
\end{proof}

\begin{prop} \label{prop-ss-quotient2}
La réduction modulo $l$ de $[\overleftarrow{t-1},\overrightarrow{s-t}]_\pi$ admet $\sharp \IC_{t-1}$ constituants
irréductibles indexés par les éléments $\underline{i} \in \IC_{t-1}$ de telle sorte que le constituant
$I_{\underline{i}}([\overleftarrow{t-1},\overrightarrow{s-t}]_\pi)$ est de niveau de cuspidalité $\underline{i}$
isomorphe à
$$I_{\underline{i}}([\overleftarrow{m(\varrho)\underline{i}(l)-1}]_\pi) \times I_{\underline{0}}
([\overleftarrow{t_\varrho(\underline{i})-1},\overrightarrow{s-t}]_\pi).$$
\end{prop}

\begin{proof}
Le raisonnement procède en deux temps: on montre tout d'abord que les représentations
$I_{\underline{i}}([\overleftarrow{m(\varrho)\underline{i}(l)-1}]_\pi) \times
I_{\underline{0}}([\overleftarrow{t_\varrho(\underline{i})-1},\overrightarrow{s-t}]_\pi)$, qui sont irréductibles
d'après \cite{vigneras-induced} \S V.3, sont effectivement des constituants de la réduction modulo $l$ de
$[\overleftarrow{t-1},\overrightarrow{s-t}]_\pi$. Ainsi
$$\Psi=r_l\Bigl ( [\overleftarrow{t-1},\overrightarrow{s-t}]_\pi \Bigr ) - \sum_{\underline{i} \in \IC_{t-1}}
I_{\underline{i}}([\overleftarrow{m(\varrho)\underline{i}(l)-1}]_\pi) \times
I_{\underline{0}}([\overleftarrow{t_\varrho(\underline{i})-1},\overrightarrow{s-t}]_\pi)$$ est un élément effectif du
groupe de Grothendieck dont on montre dans une deuxième étape, que son image par par tout foncteur de Jacquet est
nulle; le résultat découle alors du fait que $r_l\Bigl ( [\overleftarrow{t-1}, \overrightarrow{s-t}]_\pi \Bigr )$ ne
contient pas de cuspidales.

Pour tout $m(\varrho)\underline{i}(l) \leq t-1$, on a
$$[\overleftarrow{t-1},\overrightarrow{s-t}]_\pi \hookrightarrow
[\overleftarrow{t_\varrho(\underline{i})-1},\overrightarrow{s-t}]_{\pi \{ m(\varrho)\underline{i}(l)/2 \} } \times
[\overleftarrow{m(\varrho)\underline{i}(l)-1}]_{\pi \{ -s_\varrho(\underline{i})/2 \} }$$ de sorte que, de la
connaissance de
$J_{P_{s_\varrho(\underline{i})g,m(\varrho)\underline{i}(l)g}}([\overleftarrow{t-1},\overrightarrow{s-t}]_\pi)$, et de
l'exactitude du foncteur de Jacquet $J_{P_{s_\varrho(\underline{i})g,m(\varrho)\underline{i}(l)g}}$, on en déduit que
$$J_{P_{s_\varrho(\underline{i})g,m(\varrho)\underline{i}(l)g}}([\overleftarrow{t-1},\overrightarrow{s-t}]_\pi)
\twoheadrightarrow [\overleftarrow{t_\varrho(\underline{i})-1},\overrightarrow{s-t}]_{\pi \{
m(\varrho)\underline{i}(l)/2 \} } \otimes [\overleftarrow{m(\varrho)\underline{i}(l)-1}]_{\pi \{
-s_\varrho(\underline{i})/2 \} }.$$ Comme la réduction modulo $l$ commute aux foncteurs de Jacquet, et qu'il existe un
réseau de $[\overleftarrow{t_\varrho(\underline{i})-1},\overrightarrow{s-t}]_{\pi \{ m(\varrho)\underline{i}(l)/2 \} }
\otimes [\overleftarrow{m(\varrho)\underline{i}(l)-1}]_{\pi \{ -s_\varrho(\underline{i})/2 \} }$ dont la réduction
modulo $l$ est semi-simple, on en déduit, par réciprocité de Frobenius, que les représentations irréductibles
$\IC_0([\overleftarrow{t_\varrho(\underline{i})-1},\overrightarrow{s-t}]_\pi) \times \IC_{\underline{i}}
([\overleftarrow{m(\varrho)\underline{i}(l)-1}]_\pi)$ sont des constituant de $r_l\Bigl (
[\overleftarrow{t-1},\overrightarrow{s-t}]_\pi \Bigr )$. \footnote{On notera que $\varrho \{
m(\varrho)\underline{i}(l)/2 \} \simeq \varrho$}

Il s'agit désormais de montrer avec les notations ci-dessus, que pour tout parabolique $P$, $J_P(\Psi)$ est nul. On
raisonne par l'absurde: on considère alors un parabolique $P$, tel que $J_P(\Psi)$ contienne un constituant de niveau
de cuspidalité minimal. Par transitivité des foncteurs de Jacquet et en remarquant que ceux-ci augmente le niveau de
cuspidalité, on peut supposer que le facteur de Lévi de $P$ est de la forme $GL_{(s-t')g} \times GL_{t'g}$. Ce
constituant $\Pi_0$ est alors tel que son image par $J_P$ est à prendre parmi les constituants de la réduction modulo
$l$ de
$$[\overleftarrow{t-t'-1+a},\overrightarrow{s-t-a}]_{\pi \{ (t'-a)/2 \} }
\otimes [\overleftarrow{t'-a-1}]_{\pi \{ (t'-s-a)/2 \} } \times [\overrightarrow{a-1}]_{ \pi \{ (s-a)/2 \} }$$ pour $0
\leq a \leq \min \{ s-t,t' \}$. On remarque alors que $J_{P_{g,(s-1)g}}(\Pi_0)$ est non nulle: en effet pour $s-t-a>0$,
par réciprocité de Frobenius, $\Pi_0$ est un sous-espace d'une induite parabolique de la forme
$I_{\underline{i}}([\overleftarrow{t-t'-1+a},\overrightarrow{s-t-a}]_{\pi \{ ? \} } \times \Pi$ où $\Pi$ est une
certaine représentation irréductible de $GL_{t'g}(K)$ que l'on peut écrire sous la forme d'une induite $\Pi_1 \times
\Pi_2$ avec $\Pi_1$ qui est $\varrho$-superunipotente et $\Pi_2$ étant de $\varrho_{-1}$-niveau de cuspidalité nul.
Ainsi comme
$$I_{\underline{i}}([\overleftarrow{t-t'-1+a},\overrightarrow{s-t-a}]_\pi \simeq I_{\underline{0}}
([\overleftarrow{t"-1},\overrightarrow{s-t-a}]_\pi \times I_{\underline{i}}([\overleftarrow{s-t'-t"-1}]_\pi),$$ $\Pi_0$
s'écrit sous la forme d'une induite irréductible d'une représentation $\varrho$-superunipotente, obtenue comme un
sous-espace forcément $\varrho$-superunipotente de $I_{\underline{0}}([\overleftarrow{t"-1},\overrightarrow{s-t-a}]_\pi
\times \Pi_1$ avec une représentation de $\varrho_{-1}$-niveau de cuspidalité nul; on en déduit donc que
$J_{P_{g,(s-t'-1)g}}(\Pi_0)$ est non nulle.

Pour $a=s-t$, le même argument fonctionne dès que $s-t'-m(\varrho) \underline{i}(l)>0$. Dans les cas restant, $s-t'$
est divisible par $m(\varrho)$ de sorte que
$$r_l \Bigl ( [\overleftarrow{t'-a-1}]_{\pi \{ (t'-a-s)/2 \} } \times [\overrightarrow{a-1}]_{ \pi \{ (s-a)/2 \} }
\Bigr ) = r_l \Bigl ( [\overleftarrow{t'-a-1},\overrightarrow{a+1}]_{\pi \{ ? \} } \Bigr ) + r_l \Bigl (
[\overleftarrow{t'-a},\overrightarrow{a-1}]_{\pi \{ ? +1/2 \} } \Bigr )$$ et on se retrouve dans la situation
précédente, d'où le résultat.

Ainsi pour conclure, il suffit de vérifier que $J_{P_{g,(s-1)g}}(\Psi)$ est nulle. Un constituant irréductible de $r_l
\Bigl ( J_{P_{g,(s-1)g}}([\overleftarrow{t-1},\overrightarrow{s-t}]_\pi) \Bigr )$ est de la forme
$$ \varrho \{ (2t-1-s)/2 \} \otimes I_{\underline{i}}([\overleftarrow{t-2}]_{\pi\{ (t-s-1)/2 \} }) \times
\times [\overrightarrow{s-t-1}]_{\varrho\{ t/2 \} }$$ qui n'apparaît qu'avec multiplicité $1$. On remarque alors que
celui-ci est un constituant de $J_{P_{g,(s-1)g}} (I_{\underline{i}}([\overleftarrow{m(\varrho)\underline{i}(l)-1}]_\pi)
\times I_{\underline{0}}([\overleftarrow{t_\varrho(\underline{i})-1},\overrightarrow{s-t}]_\pi)$, d'où le résultat.
\end{proof}

\rem contrairement à ce qui se passe pour $t=s$, l'image des $\IC_{\underline{i}}([\overleftarrow{t-1},
\overrightarrow{s-t}]_\pi)$ par les foncteurs de Jacquet $J_P$, n'est pas égale à la somme des constituants de niveau
de cuspidalité $\underline{i}$: il ne les contient pas tous et il en contient d'autres de niveau supérieur comme le
montre l'exemple suivant.

\noindent \textit{Exemple}: la réduction modulo $l=3$ de $[\overleftarrow{1},\overrightarrow{1}]_\pi$ dans le cas où
$m(\varrho)=2$, est irréductible et donc $\varrho$-superunipotente; $J_{P_{2g,g}}(I_{\underline{0}}
([\overleftarrow{1},\overrightarrow{1}]_\pi)$ est dans le groupe de Grothendieck égal à
$$\Bigl ( I_{\underline{0}}([\overleftarrow{1}]_{\pi\{ -\frac{1}{2} \} }) + [\overrightarrow{1}]_{\varrho \{ -\frac{1}{2} \} }
+ I_{(1,0,\cdots)}([\overleftarrow{1}]_{\pi\{ -\frac{1}{2} \} }) \Bigr ) \otimes \varrho \{ \frac{1}{2} \}.$$ Dans ce
cas on peut même préciser les extensions: en effet d'après la première et la seconde formule d'adjonction, tout
sous-espace et tout quotient irréductible de $J_{P_{2g,g}}(I_{\underline{0}}
([\overleftarrow{1},\overrightarrow{1}]_\pi))=\Pi \{ -\frac{1}{2} \} \otimes \varrho \{ \frac{1}{2} \}$ doit être
$\varrho$-superunipotent de sorte que le graphe d'extensions $\Gamma_\Pi$ possède deux flèches, celle reliant
$[\overrightarrow{1}]_\varrho$ à $I_{(1,0,\cdots)}([\overleftarrow{1}]_\pi)$ et celle reliant
$I_{(1,0,\cdots)}([\overleftarrow{1}]_\pi)$ à $I_{\underline{0}}([\overleftarrow{1}]_\pi)$.

\subsection{Construction de réseaux d'induction}

Comme la réduction modulo $l$ de $[\overrightarrow{t-1}]_\pi$ est irréductible, il ne possède à isomorphisme près qu'un
seul réseau stable que l'on notera encore $[\overrightarrow{t-1}]_\pi$. Ainsi l'image de
$[\overleftarrow{t-1},\overrightarrow{s-t}]_\pi$ dans $RI_{\bar \Zm_l}([\overleftarrow{t-1}]_\pi)
\overrightarrow{\times} [\overrightarrow{s-t-1}]_\pi$ définit un réseau que l'on note $RI_{\bar
\Zm_l}([\overleftarrow{t-1},\overrightarrow{s-t}]_\pi)$. En considérant la surjection $[\overleftarrow{t-2}]_\pi
\overrightarrow{\times} [\overrightarrow{s-t}]_\pi \twoheadrightarrow [\overleftarrow{t-1},\overrightarrow{s-t}]_\pi$
permet aussi à partir du réseau $RI_{\bar \Zm_l} ([\overleftarrow{t-2}]_\pi)$ de définir un autre réseau de
$[\overleftarrow{t-1},\overrightarrow{s-t}]_\pi$. La proposition suivante montre que ce dernier réseau n'est autre que
$RI_{\bar \Zm_l}([\overleftarrow{t-1},\overrightarrow{s-t}]_\pi)$.

\begin{prop} \label{prop-compatible}
Soit $\pi$ une $\bar \Qm_l$-représentation entière cuspidale irréductible de $GL_g(K)$ et
$$f: [\overleftarrow{t-1}]_{\pi} \overrightarrow{\times} [\overrightarrow{a}]_{\pi}
\longrightarrow [\overleftarrow{t}]_{\pi} \overrightarrow{\times} [\overrightarrow{a-1}]_{\pi}$$ un morphisme
$GL_{(t+a+1)g}(K)$-équivariant non nul et défini sur $\bar \Zm_l$. Si les réseaux associés aux représentations de
Steinberg généralisées sont $RI_{\bar \Zm_l}([\overleftarrow{t-1}]_\pi)$ et $RI_{\bar \Zm_l} ([\overleftarrow{t}]_\pi)$
alors $f$ est donné par un scalaire $\lambda_f \in \bar \Zm_l$ de sorte que la réduction modulo $l$ de $f$ est soit
nulle soit d'image la réduction modulo $l$ de $[\overleftarrow{t},\overrightarrow{a}]_\pi$.
\end{prop}

\begin{proof}
Autrement dit il s'agit de montrer que les réseaux $R_{\bar \Zm_l,\pm}([\overleftarrow{t},\overrightarrow{a}]_\pi)$
induits respectivement par $RI_{\bar \Zm_l}([\overleftarrow{t-1}]_{\pi}) \overrightarrow{\times}
[\overrightarrow{a}]_{\pi}$ et $RI_{\bar \Zm_l}([\overleftarrow{t}]_{\pi}) \overrightarrow{\times}
[\overrightarrow{a-1}]_{\pi}$ sont les mêmes. D'après \ref{prop-reseau}, on a sur une $\bar \Zm_l$-surjection
$$RI_{\bar \Zm_l}([\overleftarrow{t-1}]_{\pi}) \overrightarrow{\times} [\overrightarrow{0}]_\pi
\overrightarrow{\times} [\overrightarrow{a-1}]_{\pi} \twoheadrightarrow RI_{\bar \Zm_l}([\overleftarrow{t}]_{\pi})
\overrightarrow{\times} [\overrightarrow{a-1}]_{\pi}$$ qui induit donc une $\bar \Zm_l$-surjection
$$RI_{\bar \Zm_l}([\overleftarrow{t-1}]_{\pi}) \overrightarrow{\times} [\overrightarrow{a}]_{\pi}
\twoheadrightarrow R_{\bar \Zm_l,-}([\overleftarrow{t},\overrightarrow{a}]_{\pi})$$ et donc un isomorphisme $R_{\bar
\Zm_l,+}([\overleftarrow{t},\overrightarrow{a}]_\pi) \simeq R_{\bar
\Zm_l,-}([\overleftarrow{t},\overrightarrow{a}]_\pi)$.
\end{proof}

La proposition suivante caractérise les réseaux $RI_{\bar \Zm_l}([\overleftarrow{t-1},\overrightarrow{s-t}]_\pi)$ en
précisant le graphe des extensions de sa réduction modulo $l$.

\begin{prop} \label{prop-reseaubis}
Les flèches du graphe des extensions de la réduction modulo $l$ de $RI_{\bar
\Zm_l}([\overleftarrow{t-1},\overrightarrow{s-t}]_\pi)$ sont celles qui relient deux constituants irréductibles
$I_{\underline{i}}([\overleftarrow{t-1},\overrightarrow{s-t}]_\pi)$ et
$I_{\underline{i'}}([\overleftarrow{t-1},\overrightarrow{s-t}]_\pi)$ où $\underline{i}<\underline{i'}$ sont deux
éléments consécutifs de $\IC_{t-1}$.
\end{prop}

\begin{proof}
Notons que quelque soit le constituant irréductible de $r_l \Bigl ( [\overleftarrow{t-1},\overrightarrow{s-t}]_\pi
\Bigr )$, son image par $J_{P_{(s-t+1)g,(t-1)g}}$ admet un constituant de la forme
$I_{\underline{i}}([\overleftarrow{t-1}]_{\pi \{ (t-s)/2 \} } ) \otimes [\overrightarrow{s-t-1}]_{\varrho \{ t/2 \} }$.
Par ailleurs en raisonnant comme dans la preuve de la proposition \ref{prop-reseau}, on obtient, cf. \ref{eq-choix},
une surjection
$$\phi:J_{P_{tg,(s-t)g}}(RI_{\bar \Zm_l}([\overleftarrow{t-1},\overrightarrow{s-t}]_\pi)
\twoheadrightarrow RI_{\bar \Zm_l} ([\overleftarrow{t-1}]_{\pi \{ (t-s)/2 \} } \otimes
[\overrightarrow{s-t-1}]_{\varrho \{ t/2 \} }$$ de sorte que si
$$(0)=V_n \subset V_{n-1} \subset \cdots \subset V_1 \subset V_0=r_l \Bigl ( RI_{\bar \Zm_l}([\overleftarrow{t-1},\overrightarrow{s-t}]_\pi \Bigr )$$
est une filtration dont les gradués sont irréductibles, alors
\begin{multline*}
(0)=\phi \circ J_{P_{tg,(s-t)g}}(V_n) \subset \phi \circ J_{P_{tg,(s-t)g}}(V_{n-1}) \subset \cdots \subset  \\ \phi
\circ J_{P_{tg,(s-t)g}}(V_1) \subset r_l(RI_{\bar \Zm_l}([\overleftarrow{t-1}]_{\pi \{ (t-s)/2 \} } \otimes
[\overrightarrow{s-t-1}]_{\varrho \{ t/2 \} }
\end{multline*}
est encore une filtration dont tous les gradués sont non nuls de sorte que, comme $\sharp \IC_{t-1}$ est égal au nombre
des constituants irréductibles de $r_l([\overleftarrow{t-1},\overrightarrow{s-t}]_\pi)$ et de
$r_l([\overleftarrow{t-1}]_\pi)$, tous ces gradués sont encore irréductibles. Le résultat découle alors directement de
la même propriété pour $RI_{\bar \Zm_l}([\overleftarrow{t-1}]_\pi)$.
\end{proof}

\rem $RI_{\bar \Zm_l}([\overleftarrow{t-1},\overrightarrow{s-t}]_\pi)$ est aussi le réseau induit par
$[\overrightarrow{s-t}]_\pi \overleftarrow{\times} RI_{\bar \Zm_l}([\overleftarrow{t-2}]_\pi)$ et donc d'après le lemme
\ref{lem-reseau0} de $[\overrightarrow{s-t-1}]_\pi \overleftarrow{\times} RI_{\bar \Zm_l}([\overleftarrow{t-1}]_\pi)$.
Pour le vérifier, il suffit de montrer que le graphe des extensions de la réduction modulo $l$ est celui de $RI_{\bar
\Zm_l}([\overleftarrow{t-1},\overrightarrow{s-t}]_\pi$: pour cela on raisonne comme dans la preuve ci-dessus en
utilisant $J_{P_{(s-t+1)g,(t-1)g}}$ composé avec la surjection sur $[\overrightarrow{s-t}]_{\varrho \{ (t-1)/2 \} }
\otimes RI_{\bar \Zm_l} ([\overleftarrow{t-2}]_{\pi \{ (t-s-1)/2 \} })$.

\begin{coro} Soit $K^\bullet$ le complexe tel que pour $1 \leq i \leq s$,
$$K^{d-i}=RI_{\bar \Zm_l}([\overleftarrow{s-i}]_\pi) \overrightarrow{\times} [\overrightarrow{i-2}]_\varrho$$
et $K^i$ est nul pour les autres indices. On suppose que les flèches $d_i$ pour $1 \leq i \leq s$ sont non nulles,
entières données, d'après le lemme de Schur, par un scalaire $\lambda_i \in \bar \Zm_l$ que l'on suppose en outre
inversible. En notant $U^i=\ker d_i/\im d_{i-1}$, on a alors
$$U_i=RI_{\bar \Zm_l} \Bigl ( [\overleftarrow{s-i},\overrightarrow{i-1}]_\pi).$$
\end{coro}

Le résultat précédent est motivé par le complexe de Deligne-Carayol associé aux cycles évanescents de l'espace
des déformations d'un $\OC_K$-module formel de hauteur $d=sg$: dans un récent papier, nous montrons que les
hypothèses du corollaire précédent sont vérifiées.

\subsection{La réduction modulo $l$ des représentations $LT_t(\pi,s)$ sont étrangères}

L'objet de ce paragraphe est de montrer que pour $t \neq t_1$, les réductions modulo $l$ de
$[\overleftarrow{t-1},\overrightarrow{s-t}]_\pi$ et $[\overleftarrow{t_1-1},\overrightarrow{s-t_1}]_\pi$ sont
disjointes.

\begin{lemm} \label{lem-Ztrivial2}
Si $m_l(\varrho) \neq 2$, alors pour tout $s \geq 2$ et pour tout $\delta$, $[\overrightarrow{s-1}]_{\varrho \{ \delta
\} }$ n'est pas un sous-quotient de la réduction modulo $l$ de $[\overleftarrow{s-1}]_{\pi}$.
\end{lemm}

\begin{proof} D'après (\ref{lem-Ztrivial}), si $\epsilon_l(\varrho)>2$, $I_{\underline{0}}([\overleftarrow{s-1}]_{\pi})$ n'est pas de
la forme $[\overrightarrow{s-1}]_{\varrho \{ \delta \} }$. Considérons alors le cas $\epsilon_l(\varrho)=1$ de sorte
que $m_l(\varrho)=l \neq 2$ et raisonnons par l'absurde. Par application du foncteur de Jacquet $J_{P_{2g,(s-2)g}}$, on
est ramené au cas $s=2$; la contradiction découle alors du fait que pour $l \neq 2$, la réduction modulo $l$ de
$[\overleftarrow{1}]_{\pi}$ est irréductible isomorphe à $[\overleftarrow{1}]_{\varrho}$.
\end{proof}

\begin{prop} \label{prop-triviale-non}
Soient $l \neq 2$ et $0 < t < a$. Alors pour tout $\delta$, $[\overrightarrow{a}]_{\varrho \{ \delta \} }$ n'est jamais
isomorphe à $I_{\underline{0}}([\overleftarrow{t},\overrightarrow{a-t}]_{\pi})$.
\end{prop}

\begin{proof} On raisonne par l'absurde. D'après la proposition \ref{prop-super-unique2}, par application de
$J_{P_{(t+1)g,(a-t)g}}$, on devrait alors avoir
$$[\overrightarrow{t}]_{\varrho \{ \delta+\frac{t-a}{2} \} } \simeq
I_{\underline{0}}([\overleftarrow{t}]_{\pi\{ \frac{t-a}{2} \} }) \qquad [\overrightarrow{a-t-1}]_{\varrho \{
\delta+\frac{t+1}{2}\} } \simeq [\overrightarrow{a-t-1}]_{\varrho\{ \frac{t+1}{2}\} })$$ Du deuxième isomorphisme, on
obtient $\delta \equiv 0 \mod \epsilon_l(\pi)$. On en déduit alors que $[\overrightarrow{t}]_{\varrho} \simeq
I_{\underline{0}}([\overleftarrow{t}]_{\pi})$ de sorte que d'après le lemme précédent on a $\epsilon_l(\varrho)=2$ (si
$\epsilon_l(\varrho)=1$ alors $m_l(\varrho)=l$ est distinct de $2$ par hypothèse). Par ailleurs en appliquant
$J_{P_{g,tg}}$ au premier isomorphisme, l'égalité des supports cuspidaux impose alors $t \equiv 0 \mod
\epsilon_l(\varrho)$. En considérant $J_{P_{tg,(a-t+1)g}^{op}}$, d'après (\ref{prop-super-unique2}), on obtient selon
le même procédé $a \equiv 0 \mod \epsilon_l(\varrho)$. Ainsi avec $\epsilon_l(\varrho)=2$, on a $t\geq 2$ et $a-t \geq
2$. En appliquant le foncteur de Jacquet $J_{P_{ag,g}}$, on obtient alors que $[\overrightarrow{a-1}]_{\varrho \{ -1/2
\} }$ est égal soit à $I_{\underline{0}}([\overleftarrow{t},\overrightarrow{a-t-1}]_{\pi\{ -1/2\} }$, soit à
$I_{\underline{0}}([\overleftarrow{t-1},\overrightarrow{a-t}]_{\pi\{ 1/2 \} }$. La contradiction découle alors d'un
raisonnement par récurrence.
\end{proof}

\rem en appliquant l'involution de Zelevinski et la dualité à la proposition précédente, on en déduit que
$$I_{\underline{0}}([\overleftarrow{a}]_{\pi^\vee \{ \delta \} })=\Bigl ( Z_\varrho ([\overrightarrow{a}]_{\varrho \{ - \delta \} }) \Bigr
)^\vee$$ n'est jamais isomorphe à
$$I_{\underline{0}}([\overleftarrow{a-t},\overrightarrow{t}]_\pi^\vee)=I_{\underline{0}}( Z_\pi ( [\overleftarrow{t},\overrightarrow{a-t}]_\pi))^\vee$$
pour tout $\delta$ dès que $0<t<a$.

\begin{coro} \label{coro-etrangeres}
Soient $l \neq 2$ et $0 < t<t_1 < a$; la réduction modulo $l$ de $[\overleftarrow{t-1},\overrightarrow{s-t}]_\pi$ et
$[\overleftarrow{t_1-1},\overrightarrow{s-t_1}]_\pi$ sont disjointes.
\end{coro}

\begin{proof}
Considérons pour $\underline{i} \in \IC_{t-1}$, les constituants de niveau de cuspidalité $\underline{i}$ de
$[\overleftarrow{t-1},\overrightarrow{s-t}]_\pi$ et $[\overleftarrow{t_1-1},\overrightarrow{s-t_1}]_\pi$. D'après la
proposition \ref{prop-super-unique2}, $J_{P_{(s-t+1)g,(t-1)g}}\Bigl (
I_{\underline{i}}([\overleftarrow{t-1},\overrightarrow{s-t}]_\pi) \Bigr )$ contient $[\overrightarrow{s-t}]_{\varrho \{
(t-1)/2 \} } \otimes I_{\underline{i}}([\overleftarrow{t-2}]_{\pi \{ (t-1-s)/2 \} })$ alors que les constituants de
$$J_{P_{(s-t+1)g,(t-1)g}}\Bigl ( I_{\underline{i}}([\overleftarrow{t_1-1},\overrightarrow{s-t_1}]_\pi) \Bigr )$$
tels que
le premier terme soit $\varrho$-superunipotent, sont de la forme
$I_{\underline{0}}([\overleftarrow{a},\overrightarrow{s-t-a}]_{\pi \{ \delta \} } \otimes ?$ pour $2\delta=2t_1-2a-t-1$
et avec $t_1-t \leq a \leq s-t$; d'après la proposition précédente
$I_{\underline{0}}([\overleftarrow{a},\overrightarrow{s-t-a}]_{\pi \{ \delta \} }$ ne peut pas être isomorphe à
$[\overrightarrow{s-t}]_{\varrho \{ (t-1)/2 \} }$ ce qui prouve le résultat.
\end{proof}

\bibliographystyle{plain}
\bibliography{bib-ok}

\end{document}